\documentclass[12pt,a4paper,twoside]{article}

\usepackage[active]{srcltx}
\usepackage[utf8]{inputenc}
\usepackage{hyperref}
\usepackage[english]{babel}
\usepackage{authblk}

\usepackage[a4paper,hcentering,vcentering]{geometry}
\usepackage{url,epsfig,amssymb,amsthm,amsmath}

\numberwithin{equation}{section}

\numberwithin{theorems}{section}

\numberwithin{corollarys}{section}

\newtheorem{definition}{Definition}
\numberwithin{definition}{section}

\newtheorem{lemma}{Lemma}
\numberwithin{lemma}{section}
\newtheorem{proposition}{Proposition}
\numberwithin{proposition}{section}
\newtheorem{theoremx}{Theorem}

\newtheorem{corollaryx}{Corollary}

\newtheorem*{definitionnull}{Definition}

\newcommand{\N}{\ensuremath{\mathbb{N}}}

\newcommand{\R}{\ensuremath{\mathbb{R}}}
\newcommand{\T}{\ensuremath{\mathbb{T}}}
\newcommand{\Z}{\ensuremath{\mathbb{Z}}}

\begin{document}

\title{\LARGE{\textbf{Asymptotic motions converging to arbitrary dynamics for time-dependent Hamiltonians}}}%
\author{Donato Scarcella}
\affil{Université Paris-Dauphine - Ceremade UMR 7534 Place du Maréchal De Lattre De Tassigny, 75016 PARIS.}
\date{}%
\maketitle


\begin{abstract}
In a previous work~\cite{Sca22a}, we consider time-dependent perturbations of a Hamiltonian having an invariant torus supporting quasiperiodic solutions. Assuming the perturbation decays polynomially fast as time tends to infinity, we prove the existence of an asymptotic KAM torus. That is a time-dependent family of embedded tori converging as time tends to infinity to the quasiperiodic invariant torus of the unperturbed system. 

In this paper, the dynamic on the invariant torus associated with the unperturbed Hamiltonian is arbitrary. Therefore, we need to assume exponential decay in time in order to prove the existence of a time-dependent family of embedded tori converging in time to the invariant torus associated with the unperturbed system. The proof relies on the implicit function theorem, and the most complicated and original part rests on the solution of the associated linearized problem.
\end{abstract}


\setcounter{tocdepth}{1}
\tableofcontents

\section{Introduction}

The purpose of this paper, together with~\cite{Sca22a, Sca22b, Sca22c}, is to develop a non-autonomous KAM theory. The classical KAM theory deals with the persistence of quasiperiodic solutions in nearly integrable Hamiltonian systems. The pioneering works are those of Kolmogorov~\cite{kolmogorov1954conservation}, Arnold~\cite{Arn63a, Arn63b} and Moser~\cite{moser1962invariant}. 

Over the years, these works inspired many mathematicians, enriching KAM theory with many variations and improvements. We refer to Bost~\cite{bost1986tores}, Pöschel ~\cite{Posc01}, Chierchia~\cite{Chi03} and Féjoz~\cite{Fe16} for very interesting surveys and Dumas's book~\cite{Dumas} for a fascinating historical report. 
On the other hand, a rich and exhaustive survey of this theory differs from the purpose of this introduction. For this reason, we will limit ourselves to suggesting some important works whose techniques have deeply influenced the present paper.

Pöschel~\cite{Posc01} provides a refined statement and an elegant proof of the KAM theorem in the case of real-analytic Hamiltonians. It is based on an idea of Moser~\cite{Mos67} introducing the frequencies as independent parameters. Letting an integer $n \ge 2$ and a real number $\tau > n-1$, we consider the case of finite differentiable Hamiltonians. We know that it suffices the Hamiltonian to be of class $C^k$ with $k>2(\tau+1)>2n$. Here, we refer to the works of~\cite{Posc82},~\cite{Sal04},~\cite{Bou20} and~\cite{Kou20}.
On the other hand, many interesting proofs are given using a Nash-Moser approach. It consists of the introduction of an adapted implicit function theorem in a scale of Banach spaces (replacing the iterative scheme introduced by Kolmogorov). We refer to the works of Zehnder~\cite{Zeh76, Zeh75}, Herman~\cite{bost1986tores}, Berti-Bolle~\cite{BB15} and Féjoz~\cite{Fe04}.  Concerning the dissipative case, one can see Massetti~\cite{M19}.

In this paper, we are interested in time-dependent perturbations of Hamiltonians having an invariant torus. The first papers in this direction are those of Fortunati-Wiggins~\cite{FW14} and Canadell-de la Llave~\cite{CdlL15}. The work of Canadell-de la Llave generalizes the one of Fortunati-Wiggins. They prove the existence of an asymptotic KAM torus for time-dependent vector fields converging exponentially fast as time tends to infinity to autonomous vector fields having an invariant torus supporting quasiperiodic solutions. For an asymptotic KAM torus, we mean a time-dependent family of embedded tori converging as time tends to infinity to the quasiperiodic invariant torus associated with the autonomous system. 

In another work~\cite{Sca22a}, we generalise the results of Fortunati-Wiggins and Canadell-de la Llave in the particular case of Hamiltonian systems. We consider time-dependent perturbation of real-analytic or finite differentiable Hamiltonians having an invariant torus supporting quasiperiodic solutions. We assume the perturbation decays polynomially fast as time tends to infinity without any smallness assumption on the perturbation, and we prove the existence of an asymptotic KAM torus. 

We also studied when we have the existence of a biasymptotic KAM torus. That is a continuous time-dependent family of embedded tori converging to an invariant torus in the future (when $t \to +\infty$) and an invariant torus in the past (when $t \to -\infty$). But, to be more precise, we went one step further. Using perturbation techniques, considering time-dependent perturbations of integrable Hamiltonians or Hamiltonians having a large (in the sense of measure) set of invariant tori, he proves the existence of orbits converging to quasiperiodic solutions in the future and the past~\cite{Sca22c}. 

These kinds of perturbations are not artificial and appear in many physical problems. We refer to~\cite{kawai2007transition, BdlL11} for the example of a molecule disturbed by another molecule or by a laser pulse. On the other hand, we also analysed the example of the planar three-body problem perturbed by a given comet coming from and going back to infinity asymptotically along a hyperbolic Keplerian orbit modelled as a time-dependent perturbation~\cite{Sca22b}. 

As mentioned before, in the present paper, we are interested in time-dependent perturbations of Hamiltonians having an invariant torus. In this case, we assume that the dynamic on the invariant torus associated with the unperturbed system is arbitrary (and hence not quasiperiodic). In this case, similarly to Fortunati-Wiggins and Canadell-de la Llave, we assume that the perturbation decays exponentially fast in time, and we prove the existence of an asymptotic torus. That is a time-dependent family of embedded tori converging as time tends to infinity to the invariant torus associated with the unperturbed system that supports an arbitrary dynamic. 

For the sake of clarity, let us introduce the definition of $C^\sigma$-asymptotic torus. Let $B \subset \R^n$ be a ball centred at the origin, and $\mathcal{P}$ is equal to $\T^n$ and $\T^n \times B$. Given $\sigma \ge 0$, $\upsilon \ge 0$ and a positive integer $k \ge 0$, we consider time-dependent vector fields $X^t$ and $X^t_0$ of class $C^{\sigma + k}$ on $\mathcal{P}$, for all $t \in J_\upsilon$, an embedding $\varphi_0$ from $\T^n $ to $\mathcal{P}$ of class $C^{\sigma}$ and a vector field on the torus $W$ of class $C^\sigma$ such that
\begin{align}
\label{hyp1AT}
& \displaystyle \lim_{t \to +\infty}  |X^t - X^t_0|_{C^{\sigma +k}} = 0,\\
\label{hyp2AT}
& X (\varphi_0(q),t) = \partial_q \varphi_0(q) W(q) \hspace{2mm} \mbox{for all $(q, t) \in \T^n \times J_\upsilon$.}
\end{align}
\begin{definition}
\label{asymKAMNC}
We assume that $(X, X_0, \varphi_0, W)$ satisfy~\eqref{hypNC1AKAM} and~\eqref{hypNC2AKAM}.
A family of $C^\sigma$ embeddings $\varphi^t: \T^n \to \mathcal{P}$ is a $C^\sigma$-asymptotic torus associated  to $(X, X_0, \varphi_0, W)$ if there exists $\upsilon' \ge \upsilon \ge 0$ such that 
\begin{align}
\label{hyp2NC}
&  \lim_{t \to +\infty}  |\varphi^t - \varphi_0|_{C^\sigma} = 0,\\
\label{hyp1NC}
&   X(\varphi(q, t), t) =  \partial_q \varphi(q, t)W(q)  + \partial_t \varphi(q, t),
\end{align}
for all $(q, t) \in \T^n \times J_{\upsilon'}$. When $\mathrm{dim} \mathcal{P} = 2n$, then $\varphi^t$ is Lagrangian if $\varphi^t(\T^n)$ is Lagrangian for all $t$.
\end{definition}

First, we want to emphasize that taking $W(q) \equiv W \in \R^n$ constant, we obtain the definition of $C^\sigma$-asymptotic KAM torus. Moreover, in the previous definition, $\mathcal{P}$ is equal to $\T^n$ and $\T^n \times B$ because we prove results for Hamiltonian systems or vector fields on the torus. 

Contrary to~\cite{Sca22a}, in this paper, we assume that the orbits associated with the unperturbed system are arbitrary. On the other hand, we need to assume that the perturbation decays exponentially fast in time (in the other work, we required polynomial decay). Similarly to~\cite{Sca22a}, the proof relies on the implicit function theorem. The most original and complicated part consists in solving and estimating the associated linearized problem (see Section \ref{SecHENC}). It is solved thanks to a suitable change of coordinate which rectifies the dynamics on the torus. The estimation requires tame estimates for the product and the composition of Hölder functions. We will prove the existence of $C^\sigma$-asymptotic tori for time-dependent Hamiltonian (Theorem \ref{Thm1NC}) and time-dependent vector field on the torus (Corollary \ref{CorNC}).

\section{Results}\label{ResultNC}

In order to state the main results of this paper, we need to introduce some notations and definitions. We recall that $B \subset \R^n$ is a ball around the origin. Moreover, for a real parameter $\upsilon \ge 0$, we define the following interval $J_\upsilon =[\upsilon, +\infty) \subset \R$. 

For each function $f$ defined on $\T^n \times B \times J_\upsilon$ and for fixed $t \in J_\upsilon$, we denote by $f^t$ the function defined on $\T^n \times B$ such that 
\begin{equation*}
f^t(q,p) = f(q,p,t)
\end{equation*}
for all $(q,p) \in \T^n \times B$. On the other hand, for fixed $p \in B$, we consider $f_{p}$ as the function defined on $\T^n \times J_\upsilon$ such that 
\begin{equation*}
f_{p}(q,t) = f(q,p,t)
\end{equation*}
for all $(q,t) \in \T^n \times J_\upsilon$. As one can expect, for fixed $(p,t) \in B \times J_\upsilon$, let $f_p^t$ be the function defined on $\T^n$ such that
\begin{equation*}
f^t_{p}(q) = f(q,p,t)
\end{equation*}
for all $q \in \T^n$. We will use this notation for the rest of this work. 

Now, given a positive real parameters $\sigma \ge 0$ and $\upsilon \ge 0$, we define the following space of functions
\begin{definition}
\label{S}
Let $\mathcal{S}^\upsilon_{\sigma}$ be the space of functions $f$ defined on  $\T^n \times B \times J_\upsilon$ such that $f \in C( \T^n \times B \times J_\upsilon)$  and, for all $t \in J_\upsilon$, $f^t \in C^\sigma( \T^n \times B)$.
\end{definition}
We use this notation also for functions defined on $\T^n \times J_\upsilon$. This will be specified by the context. For every $f \in \mathcal{S}^\upsilon_\sigma$ and for fixed $\lambda \ge 0$, we introduce the following norm
\begin{equation}
\label{normNC}
|f|^\upsilon_{\sigma,\lambda} = \sup_{t \in J_\upsilon}|f^t|_{C^\sigma}e^{\lambda t}.
\end{equation}

We point out that Section \ref{FS} is dedicated to a series of properties concerning the previous norm. 
At this moment, we need to define another space of functions. Given $\sigma$, $\upsilon \ge 0$ and an integer $k \ge 0$, we have the following definition
\begin{definition}
\label{S2}
Let $\mathcal{\bar S}^\upsilon_{\sigma, k}$ be the space of functions $f$ such that 
\begin{equation*}
f \in \mathcal{S}^\upsilon_{\sigma +k} \hspace{2mm} \mbox{and} \hspace{2mm} \partial^i_{(q, p)}f \in \mathcal{S}^\upsilon_{\sigma +k-i}
\end{equation*}
for all $0 \le i \le k$. 
\end{definition}
Here, $\partial^i_{(q, p)}$ stands for partial derivatives of order $i$ with respect to the variables $(q, p)$. Moreover, as a convention, we consider $f = \partial^0_{(q, p)}f$. In other words, the space $\mathcal{\bar S}^\upsilon_{\sigma, k}$ is composed of the functions $f \in \mathcal{S}^\upsilon_{\sigma +k}$ with partial derivatives with respect to $(q,p)$ continuous until the order $k$. 

Now, we consider a vector field on the torus $W$. We define $\mathcal{K}_W$ as the set of the Hamiltonians  $h:\T^n \times B \times J_0 \to \R$ such that, for all $(q,t) \in \T^n \times J_0$,
\begin{equation*}
h(q,0,t) = c, \quad \partial_p h(q,0,t) = W(q)
\end{equation*}
for some $c \in \R$. It is obvious that, for all $h \in \mathcal{K}_W$, the trivial embedding $\varphi_0$ given by
\begin{equation*}
\varphi_0 : \T^n \to \T^n \times B, \quad \varphi_0(q) = (q,0),
\end{equation*}
is an invariant torus for $X_h$ and the restricted vector field is $W$.

Finally, we have everything we need to state the main result of this paper.
Given $\sigma\ge 1$ and $\lambda \ge 0$, let $H$ be the Hamiltonian of the following form
\begin{equation}
\label{H1NC}
\begin{cases}
H:\T^n \times B \times J_0 \to \R\\
H(q,p,t) = h(q,p,t) + f(q,p,t)\\
W \in C^{\sigma+2}(\T^n), \quad h \in \mathcal{K}_W,\\
f_0, \partial_p f_0, \partial^2_p H \in \mathcal{\bar S}^0_{\sigma, 2},\\
|f_0|^0_{\sigma + 2, 0} + |\partial_q f_0|^0_{\sigma + 1, \lambda} < \infty, \quad |\partial_p f_0|^0_{\sigma + 2, \lambda}< \infty\\
 |\partial_p^2 H|^0_{\sigma + 2,0} < \infty
\end{cases}
\tag{$*$}
\end{equation}
\begin{theoremx}
\label{Thm1NC}
Let $H$ be as in~\eqref{H1NC}. Then, there exists a Hamiltonian $\tilde h \in \mathcal{K}_W$ and a constant $C(\sigma)$ depending on $\sigma$ such that if 
\begin{equation}
\label{lambdaNC}
\lambda > C(\sigma) |\partial_qW|_{C^0},
\tag{$\#$}
\end{equation}
there exists a Lagrangian $C^\sigma$-asymptotic torus associated to $(X_H, X_{\tilde h}, \varphi_0, W)$.
\end{theoremx}

We want to point out that the Hamiltonian $\tilde h$ is introduced in~\eqref{tildeh}. Concerning the proof, it rests on the implicit function theorem. The most complicated and original part relies on the solution of the associated linearized problem (homological equation), where~\eqref{lambdaNC} plays a crucial role (see Section \ref{SecHENC}). Indeed, the previous constant $C(\sigma)$ is specified in Lemma \ref{homoeqlemmaNC}.

We want to emphasize that our proof does not work for $C^\infty$ or real analytic Hamiltonians. The point is that we are not able to find $C^\infty$ or holomorphic solutions to the associated linearized problem (we refer to Section \ref{SecHENC} for more details).  

Contrary to the theorem of Canadell-de la Llave~\cite{CdlL15}, the vector field $W$ is not constant. In addition, by letting $W\equiv \mathrm{cst}$, we obtain $\lambda >0$, which is the hypothesis of  Canadell-de la Llave in the case of Hamiltonian systems. On the other hand, we proved~\cite{Sca22a} that we do not need exponential decay in this case. 

Here, we show the existence of a $C^\sigma$-asymptotic torus $\varphi^t$ of the form
\begin{equation*}
\varphi^t(q) = (q + u^t(q), v^t(q))
\end{equation*}
for all $q \in \T^n$ and $t$ sufficiently large, where $\mathrm{id} + u^t$ is a diffeomorphism of the torus for all fixed $t$.  Moreover, we give some information concerning the time decay of $u$ and $v$. We prove that
\begin{equation*}
|u^t|_{C^\sigma} \le C e^{-\lambda t}, \quad |v^t|_{C^\sigma} \le C e^{-\lambda t},
\end{equation*}
for all $t$ large enough and for a suitable constant $C$. 

In the case of time-dependent perturbations of vector fields on the torus, we have the following result. Given $\sigma \ge 1$, let $Z$ be a non-autonomous vector field on $\T^n \times J_0$ of the form
\begin{equation*}
\label{ZNC}
\begin{cases}
Z(q,t) = W(q) + P(q,t)\\
W \in C^{\sigma +1}(\T^n), \quad P \in \mathcal{\bar S}^0_{\sigma, 1},\\
|P|^0_{\sigma+1, \lambda} < \infty.
\end{cases}
\tag{$Z$}
\end{equation*}
\begin{corollaryx}
\label{CorNC}
Let $Z$ be as in~\eqref{ZNC}. Then, there exists a constant $C(\sigma)$ depending on $\sigma$ such that if 
\begin{equation*}
\lambda > C(\sigma) |\partial_qW|_{C^0},
\end{equation*}
there exists a $C^\sigma$-asymptotic torus $\psi^t$ associated to $(Z, W, \mathrm{Id}, W)$. 
\end{corollaryx}

\section{Functional setting}\label{FS}

This section is devoted to some properties of the norms~\eqref{normNC} introduced in the previous section. We begin by recalling the definition of the Hölder classes of functions $C^\sigma$ and some properties. 

Let $D$ be an open subset of $\R^n$ and $k \ge 0$ a positive integer. We define $C^k(D)$ as the spaces of functions $f: D \to \R$ with continuous partial derivatives $\partial^\alpha f \in C^0(D)$ for all $\alpha \in \N^n$ with $|\alpha|=\alpha_1+...+\alpha_n \le k$. Moreover, for all $f \in C^k(D)$, we have the following norm
\begin{equation*}
|f|_{C^k} = \sup_{|\alpha|\le k}|\partial^\alpha f|_{C^0},
\end{equation*}
where $|\partial^\alpha f|_{C^0} = \sup_{x \in D}|\partial^\alpha f(x)|$ denotes the sup norm. Given $\sigma=k+\mu$, with $k \in \Z$, $k \ge 0$ and $0 < \mu <1$, we define the Hölder spaces $C^\sigma(D)$ as the spaces of functions $f\in C^k(D)$ verifying
\begin{equation}
\label{Holdernorm}
|f|_{C^\sigma} = \sup_{|\alpha|\le k}|\partial^\alpha f|_{C^0} + \sup_{|\alpha| = k}{|\partial^\alpha f(x) - \partial^\alpha f(y)| \over |x-y|^\mu}<\infty.
\end{equation}
The following proposition recalls a series of well-known properties. In this paper, we denote by $C(\cdot)$ constants depending on $n$ and the other parameters in brackets.  On the other hand, $C$ denotes constants depending only on $n$.

\begin{proposition}
\label{Holder}
We consider $f$, $g \in C^\sigma(D)$ and $\sigma \ge 0$.
\begin{enumerate}
\item For all $\beta \in \N^{n}$, if $|\beta| + s = \sigma$ then  $\left|{\partial^{|\beta|} \over \partial{x_1}^{\beta_1}... \partial{x_n}^{\beta_n}} f \right|_{C^s} \le |f|_{C^\sigma}$.\\
\item  $|fg|_{C^\sigma} \le C(\sigma)\left(|f|_{C^0}|g|_{C^\sigma} + |f|_{C^\sigma}|g|_{C^0}\right)$. 
\end{enumerate}
Concerning composite functions. Let $z$ be defined on $D_1 \subset \R^n$ and takes its values on $D_2 \subset \R^{n}$ where $f$ is defined. 
If $\sigma \ge 1$ and $f \in C^\sigma (D_2)$, $z \in C^\sigma (D_1)$ then $f\circ z \in C^\sigma(D_1)$ 
\begin{enumerate}
\item[3.] $|f \circ z|_{C^\sigma} \le C(\sigma) \left(|f|_{C^\sigma}|\nabla z|^\sigma_{C^0} + |f|_{C^1}| \nabla z|_{C^{\sigma-1}}+ |f|_{C^0}\right)$.
\end{enumerate}
\end{proposition}
\begin{proof}
We refer to~\cite{Hor76} for the proof of \textit{1.} and \textit{2.}. Concerning the last property, it is proved in~\cite{Sca22a}.
\end{proof}

The norm defined by~\eqref{normNC} satisfies the following properties. As one can expect, there is a significant similarity with those enumerated in the previous proposition. 

\begin{proposition}
\label{normpropertiesNC}
For all $f$, $g \in \mathcal{S}^\upsilon_\sigma$ and positive parameters $m$, $d \ge 1$, we have the following properties.
\begin{enumerate}
\item[a.] For all $\beta \in \N^{2n}$, if $|\beta| + r \le \sigma$ then  $\left|{\partial^{|\beta|} \over \partial{q_1}^{\beta_1}... \partial{q_n}^{\beta_n} \partial{p_1}^{\beta_{n+1}}... \partial{p_n}^{\beta_{2n}}} f \right|^\upsilon_{r, \lambda} \le |f|^\upsilon_{\sigma, \lambda}$\\
\item[b.]  $|f|^\upsilon_{\sigma, \lambda} \le |f|^\upsilon_{\sigma, k\lambda}$ \\
\item[c.]  $|f g|^\upsilon_{\sigma, d\lambda+m\lambda} \le C(\sigma)\left(|f|^\upsilon_{0,d\lambda}|g|^\upsilon_{\sigma,m\lambda} + |f|^\upsilon_{\sigma,d\lambda}|g|^\upsilon_{0,m\lambda}\right)$. 
\end{enumerate}
Given $\sigma \ge 1$, for all $f$, $z \in \mathcal{S}^\upsilon_\sigma$ then $f \circ z \in \mathcal{S}^\upsilon_\sigma$ 
\begin{enumerate}
\item[d.] $|f \circ z|^\upsilon_{\sigma, k\lambda+m\lambda} \le C(\sigma) \left(|f|^\upsilon_{\sigma,k\lambda}\left(|\nabla z|^\upsilon_{0,m\lambda}\right)^\sigma + |f|^\upsilon_{1,k\lambda}|\nabla z|^\upsilon_{\sigma-1,m\lambda} +  |f|^\upsilon_{0, k\lambda+m\lambda}  \right)$,
\end{enumerate}
\end{proposition}
\begin{proof}
The proof is a straightforward application of Proposition \ref{Holder}. Properties \textit{a.} and \textit{b.} are obvious. We verify the others.

\textit{c}.
\begin{eqnarray*}
|f g|^\upsilon_{\sigma, k\lambda+m\lambda} &=& \sup_{t \in J_\upsilon}|f^t g^t|_{C^\sigma}e^{(k \lambda + m \lambda)t} \\
&\le& C(\sigma)\sup_{t \in J_\upsilon}\left(|f^t |_{C^0}|g^t |_{C^\sigma} + |f^t |_{C^\sigma}|g^t |_{C^0}\right)e^{(k \lambda + m \lambda)t} \\
&\le& C(\sigma) \sup_{t \in J_\upsilon}\left(|f^t |_{C^0}e^{k \lambda t}|g^t |_{C^\sigma}e^{m \lambda t} + |f^t |_{C^\sigma}e^{k \lambda t}|g^t |_{C^0}e^{m \lambda t}\right)\\
&\le& C(\sigma)\left(|f|^\upsilon_{0,k\lambda}|g|^\upsilon_{\sigma,m\lambda} + |f|^\upsilon_{\sigma,k\lambda}|g|^\upsilon_{0,m\lambda}\right)
\end{eqnarray*}

\textit{d}.
\begin{eqnarray*}
|f \circ z|^\upsilon_{\sigma, k\lambda+m\lambda}&=& \sup_{t \in J_\upsilon}|f^t  \circ z^t|_{C^\sigma}e^{(k \lambda + m \lambda)t} \\
&\le& C(\sigma)\sup_{t \in J_\upsilon}\left(|f^t|_{C^\sigma}e^{k\lambda t}\left(|\nabla z^t|_{C^0}e^{ m\lambda t}\right)^\sigma  e^{(1-\sigma)m \lambda t}\right)\\
&+& C(\sigma)\sup_{t \in J_\upsilon}\left(|f^t |_{C^1}e^{k\lambda t}|\nabla z^t |_{C^{\sigma-1}} e^{m\lambda t} + |f|_{C^0}e^{(k \lambda + m \lambda)t} \right)\\
&\le& C(\sigma) \left(|f|^\upsilon_{\sigma,k\lambda}\left(|\nabla z|^\upsilon_{0,m\lambda}\right)^\sigma + |f|^\upsilon_{1,k\lambda}|\nabla z|^\upsilon_{\sigma-1,m\lambda} +  |f|^\upsilon_{0, k\lambda+m\lambda}  \right)
\end{eqnarray*}
where we observe that if $t \ge 0$ and $\sigma \ge 1$ then $e^{(1-\sigma)m \lambda t} \le 1$.
\end{proof}

\section{$C^\sigma$-asymptotic torus}\label{Csigmatorus}

Here, we recall the definition of $C^\sigma$-asymptotic torus and discuss some properties of this object. Let $B \subset \R^n$ be a ball centred at the origin, $\mathcal{P}$ be equal to $\T^n$ or $\T^n \times B$ and, for all $\upsilon \ge 0$, $J_\upsilon = [\upsilon, +\infty) \subset \R$.

Given $\sigma \ge 0$, $\upsilon \ge 0$ and a positive integer $k \ge 0$, we consider time-dependent vector fields $X^t$ and $X^t_0$ of class $C^{\sigma + k}$ on $\mathcal{P}$, for all $t \in J_\upsilon$, an embedding $\varphi_0$ from $\T^n $ to $\mathcal{P}$ of class $C^{\sigma}$ and a vector field on the torus $W$ of class $C^\sigma$ such that
\begin{align}
\label{hypNC1AKAM}
& \displaystyle \lim_{t \to +\infty}  |X^t - X^t_0|_{C^{\sigma +k}} = 0,\\
\label{hypNC2AKAM}
& X (\varphi_0(q),t) = \partial_q \varphi_0(q) W(q) \hspace{2mm} \mbox{for all $(q, t) \in \T^n \times J_\upsilon$.}
\end{align}
\begin{definitionnull}[Definition \ref{asymKAMNC}]
We assume that $(X, X_0, \varphi_0, W)$ satisfy~\eqref{hypNC1AKAM} and~\eqref{hypNC2AKAM}.
A family of $C^\sigma$ embeddings $\varphi^t: \T^n \to \mathcal{P}$ is a $C^\sigma$-asymptotic torus associated  to $(X, X_0, \varphi_0, W)$ if there exists $\upsilon' \ge \upsilon \ge 0$ such that 
\begin{align}
\label{hyp2NC}
&  \lim_{t \to +\infty}  |\varphi^t - \varphi_0|_{C^\sigma} = 0,\\
\label{hyp1NC}
&   X(\varphi(q, t), t) =  \partial_q \varphi(q, t)W(q)  + \partial_t \varphi(q, t),
\end{align}
for all $(q, t) \in \T^n \times J_{\upsilon'}$. When $\mathrm{dim} \mathcal{P} = 2n$, then $\varphi^t$ is Lagrangian if $\varphi^t(\T^n)$ is Lagrangian for all $t$.
\end{definitionnull}
 As mentioned before, it generalizes the definition of $C^\sigma$-asymptotic KAM torus due to Canadell-de la Llave~\cite{CdlL15}. As one can expect, we have a series of proprieties in common with $C^\sigma$-asymptotic KAM tori. Let $\psi_{t_0, X}^t$ and $\psi_{t_0, W}^t$ be the flow at time $t$ with initial time $t_0$ of $X$ and $W$, respectively. 
 
 \begin{proposition}
 \label{prop1asymtori}
 If the flow $\psi_{t_0, X}^t$ is defined for all $t$, $t_0 \in J_{\upsilon'}$, then~\eqref{hyp1NC} is equivalent to 
 \begin{equation}
\label{hyp1bissNC}
\psi_{t_0, X}^t \circ \varphi^{t_0}(q) = \varphi^t \circ \psi_{t_0,W}^t(q)
\end{equation}
for all $t$, $t_0 \in J_{\upsilon'}$ and $q \in \T^n$.
 \end{proposition}
 \begin{proof}
 Similarly to the proof of Proposition $3.1$ in~\cite{Sca22a}, we have the claim. 
 \end{proof}

This proposition emphasizes that the orbits on this family of embeddings converge to the arbitrary dynamics associated with the unperturbed system when time tends to infinity. 
 As one can expect, the previous proposition is the key to proving that the condition~\eqref{hyp1NC} is trivial. 
 
 \begin{proposition}
 If $\psi_{t_0, X}^t$ is defined for all $t$, $t_0 \in J_{\upsilon'}$, it is always possible to find a family of embeddings satisfying~\eqref{hyp1NC}.
 \end{proposition}
 \begin{proof}
We consider an embedding $\hat \varphi : \T^n \to \mathcal{P}$. Then, we define
\begin{equation*}
\varphi^t(q) = \psi_{t_0, X}^t \circ \hat \varphi \circ \psi^{t_0}_{t,W}(q)
\end{equation*} 
for all $t$, $t_0 \in J_{\upsilon'}$ and $q \in \T^n$. The latter is a family of embeddings satisfying~\eqref{hyp1bissNC} and hence~\eqref{hyp1NC}. 
 \end{proof}
 
In the following proposition, we will see that if we have the existence of a $C^\sigma$-asymptotic torus defined for all $t$ large, then we can extend the set of definition for all $t \in \R$. 

\begin{proposition}
We assume that $\psi^t_{t_0,X}$ is defined for all $t$, $t_0 \in \R$. If there exists a $C^\sigma$-asymptotic torus $\varphi^t$ defined for all $t \ge \upsilon'$, then we can extend the set of definition for all $t \in \R$.
\end{proposition}
\begin{proof}
For all $q \in \T^n$, we define
\begin{equation*}
\phi^t(q) = \begin{cases} \varphi^t(q) \hspace{4mm} \mbox{for all $t \ge \upsilon'$}\\
\psi_{\upsilon', X}^t \circ \varphi^{\upsilon'}\circ \psi_{t, W}^{\upsilon'}(q) \hspace{4mm} \mbox{for all $t \le \upsilon'$}.\end{cases}
\end{equation*}
This is a family of embeddings that verify~\eqref{hyp2NC} and~\eqref{hyp1NC}.
\end{proof}

\section{Proof of Theorem \ref{Thm1NC}}\label{ProofThm1NC}

The idea of the proof is the same as that in~\cite{Sca22a}, except for the solution of the homological equation (here is considerably more complicated). The proof rests on the implicit function theorem. 

First, we expand the Hamiltonian $H$ in~\eqref{H1NC} in a small neighbourhood of $0 \in B$
\begin{eqnarray*}
h(q,p,t) &=& h(q,0,t) + \partial_p h(q,0,t) \cdot p +  \int_0^1 (1 -\tau) \partial^2_p h(q, \tau p,t)d\tau \cdot p^2\\
f(q,p,t) &=& f(q, 0, t) + \partial_p f(q, 0,t) \cdot p +  \int_0^1 (1 -\tau)\partial^2_p f(q, \tau p, t) d\tau\cdot p^2.
\end{eqnarray*}
We consider $h(q,0,t)=0$ for all $(q,t) \in \T^n \times J_0$, we can do it without loss of generality. Now, we denote
\begin{eqnarray*}
a(q,t) &=& f(q,0,t)\\
b(q,t) &=& \partial_p f(q,0,t)\\
m(q,p,t) &=& \int_0^1 (1 -\tau) \left(\partial_p^2 h(q,\tau p, t) + \partial_p^2 f(q, \tau p, t)\right) d\tau \\
&=& \int_0^1 (1 -\tau) \partial_p^2 H(q,\tau p, t) d\tau,
\end{eqnarray*}
and for a suitable positive real parameter $\Upsilon \ge 1$, we can rewrite the Hamiltonian $H$ in the following form
\begin{equation}
\label{H2NC}
\begin{cases}
H : \T^n \times B \times J_0 \longrightarrow \R\\
H(q,p,t) = W(q) \cdot p + a(q,t) + b(q,t) \cdot p + m(q,p,t) \cdot p^2,\\
a, b, \partial^2_p H \in \mathcal{\bar S}^0_{\sigma,2},\quad W \in C^{\sigma+2}\\
|a|^0_{\sigma + 2, 0} + |\partial_q a|^0_{\sigma + 1, \lambda} \le \Upsilon, \quad |b|^0_{\sigma + 2, \lambda} \le \Upsilon\\
 |\partial_p^2 H|^0_{\sigma + 2,0} \le \Upsilon
\end{cases}
\tag{$**$}
\end{equation}
This Hamiltonian is our new starting point. At this point, let $\tilde h$ be the following Hamiltonian
\begin{equation}
\label{tildeh}
\tilde h(q,p,t) = W(q) \cdot p + m(q,p,t) \cdot p^2
\end{equation}
for all $(q,p,t) \in \T^n \times B \times J_0$. We can see that $\tilde h \in \mathcal{K}_W$, furthermore $X_H$ and $X_{\tilde h}$ satisfy~\eqref{hyp1NC}. 

\subsection{Outline of the proof of Theorem \ref{Thm1NC}}\label{OTPNC}

We are looking for a $C^\sigma$-asymptotic torus $\varphi^t$ associated to $(X_H, X_{\tilde h}, \varphi_0, W)$, where $H$ is the Hamiltonian in~\eqref{H2NC}, $\tilde h$ is defined by~\eqref{tildeh} and $\varphi_0$ is the trivial embedding $\varphi_0 : \T^n  \to \T^n \times B$, $\varphi_0(q) = (q,0)$ . More concretely, for given $H$, we are looking for $\upsilon' \ge 0$ large  enough and some functions $u$, $v : \T^n \times J_{\upsilon'} \to \R^n$ such that 
\begin{equation*}
\varphi(q,t) = (q + u(q,t), v(q,t)),
\end{equation*}
and to satisfy the following conditions
\begin{align}
\label{hyp1NC2}
&   X_H(\varphi(q, t), t) -  \partial_q \varphi(q, t) W(q) - \partial_t \varphi(q, t) = 0,\\
\label{hyp2NC2}
&  \lim_{t \to +\infty}  |u^t|_{C^\sigma} = 0, \quad \lim_{t \to +\infty}  |v^t|_{C^\sigma} = 0,
\end{align}
for all $(q,t) \in \T^n \times J_{\upsilon'}$. The parameter $\upsilon'$ is free, and we will fix it large enough in Lemma \ref{lemmautilethmNC}.

To this end, we introduce a suitable functional $\mathcal{F}$ given by~\eqref{hyp1NC2}. First, we define
\begin{eqnarray*}
\bar m(q,p,t) p &=& \left(\int_0^1  \partial_p^2 H(q,\tau p, t) d\tau \right) p = \partial_p \Big(m(q,p,t) \cdot p^2 \Big),
\end{eqnarray*}
for all $(q,p,t) \in \T^n \times J_{\upsilon'}$. This is well defined, we refer to~\cite{Sca22a} for more details. Moreover, we introduce
\begin{equation*}
\tilde \varphi(q,t) = (q + u(q,t), v(q,t), t), \quad \tilde u(q,t) = (q + u(q,t), t),
\end{equation*}
for all $(q,t) \in \T^n \times J_{\upsilon'}$. We observe that the composition between $X_H$ and $\tilde \varphi$ equals
\begin{equation*}
X_H \circ \tilde \varphi = \begin{pmatrix}W \circ  (\mathrm{id} + u) + b \circ \tilde u + \left(\bar m \circ \tilde \varphi \right) v  \\
-\partial_q a\circ \tilde u  - \left(\partial_q W \circ  (\mathrm{id} + u) + \partial_q b \circ \tilde u\right) v - \left(\partial_q m \circ \tilde \varphi\right) \cdot v^2\end{pmatrix}
\end{equation*}
and moreover,
\begin{equation*}
\partial_q \varphi W + \partial_t \varphi = \begin{pmatrix} W + \partial_q u W + \partial_t u  \\
\partial_q vW + \partial_t v \end{pmatrix}.
\end{equation*}
The above equations are composed of functions defined on $(q,t) \in \T^n \times J_{\upsilon'}$ or $q \in \T^n$. We have omitted the dependence of the variables $(q,t)$ or $q$ to obtain a more elegant form. We keep this notation for the rest of this work. 

Now, letting 
\begin{equation*}
 \nabla u(q,t) \bar W(q) =  \partial_q u (q,t)W(q) + \partial_t u (q,t), \quad \nabla v(q,t) \bar W(q) = \partial_q v(q,t) W(q) + \partial_t v(q,t)
\end{equation*}
for all $(q, t) \in \T^n \times J_{\upsilon'}$, we can rewrite~\eqref{hyp1NC2} in the following form
\begin{eqnarray*}
\label{FfinaleOutline}
\begin{pmatrix}W \circ  (\mathrm{id} + u) - W + b \circ \tilde u + \left(\bar m \circ \tilde \varphi  \right) v - \left(\nabla u \right)\bar W \\
-\partial_q a\circ \tilde u  - \left(\partial_q W \circ (\mathrm{id} + u) + \partial_q b \circ \tilde u \right) v - \left(\partial_q m \circ \tilde  \varphi   \right) \cdot v^2 - \left(\nabla v \right)\bar W. 
\end{pmatrix} = \begin{pmatrix} 0\\0\end{pmatrix}.
\end{eqnarray*}

Thanks to the latter, we can define the above-mentioned functional $\mathcal{F}$ on suitable Banach spaces that we will specify later. Hence, let $\mathcal{F}$ be the following functional 
\begin{equation*}
\mathcal{F}(a, b, m, \bar m, W,u, v) = (F_1(b, \bar m, W,u,v), F_2(a,b, m,W,u,v))
\end{equation*}
with 
\begin{eqnarray*}
F_1(b, \bar m,W,u,v) &=& W \circ  (\mathrm{id} + u) - W + b \circ \tilde u + \left(\bar m \circ \tilde \varphi  \right) v - \left(\nabla u \right)\bar W \\,
F_2(a,b, m,W,u,v) &=& \partial_q a\circ \tilde u  + \left(\partial_q W \circ (\mathrm{id} + u) + \partial_q b \circ \tilde u \right) v\\ 
&+& \left(\partial_q m \circ \tilde  \varphi   \right) \cdot v^2 + \left(\nabla v \right)\bar W.
\end{eqnarray*}
We observe that for all $m$, $\bar m$ and $W$, 
\begin{equation*}
\mathcal{F}(0, 0, m, \bar m,W,0, 0)=0.
\end{equation*}
Thus, we can reformulate this problem in the following terms. For fixed $m$, $\bar m$ and $W$ in suitable Banach spaces and for $(a,b)$ sufficiently close to $(0,0)$, we are looking for some functions $u$, $v$ satisfying~\eqref{hyp2NC2} such that \\$\mathcal{F}(a, b, m, \bar m,W,u, v) = 0$.

The key point of the proof of Theorem \ref{Thm1NC} concerns the analysis of the associated linearized problem. One can see that the differential of $\mathcal{F}$ with respect to the variables $(u,v)$ calculated in $(0,0,m, \bar m,W, 0,0)$ is equal to
\begin{equation*}
D_{(u,v)} \mathcal{F}(0,0,m, \bar m,W, 0,0) (\hat u, \hat v) = \begin{pmatrix} \partial_q W\hat u - \left(\nabla \hat u\right)\bar W + \bar m_0\hat v\\
\partial_q W\hat v +\left(\nabla \hat v\right)\bar W \end{pmatrix},
\end{equation*}
where, by the notation introduced at the beginning of Section \ref{ResultNC}, $\bar m_0(q,t) = \bar m(q,0,t)$ for all $(q,t) \in \T^n \times J_{\upsilon'}$. The proof of the existence of a right inverse of the latter contains the most original and mathematically complicated part of this paper.

In the following four sections, we prove Theorem \ref{Thm1NC}. First, we introduce suitable Banach spaces to properly define the functional $\mathcal{F}$. Then, we solve the homological equation, and we prove that $D_{(u,v)} \mathcal{F}(0,0,m, \bar m,W, 0,0)$ admits a right inverse. In the penultimate section, we verify that $\mathcal{F}$ satisfies the hypothesis of the implicit function theorem, and in the last section, we conclude the proof.  

\subsection{Preliminary settings}\label{PSNC}

Let $\sigma$, $\lambda$ and $\Upsilon$ be the positive parameters introduced by~\eqref{H1NC}. For $\upsilon' \ge 0$ that we will specify later, we consider the following Banach spaces $\left(\mathcal{A}, |\cdot |\right)$, $\left(\mathcal{B}, |\cdot |\right)$, $\left(\mathcal{U}, |\cdot |\right)$, $\left(\mathcal{V}, |\cdot |\right)$, $\left(\mathcal{Z}, |\cdot |\right)$ and $\left(\mathcal{G}, |\cdot |\right)$ 
\vspace{5mm}
\begin{eqnarray*}
\mathcal{A} &=& \Big\{a : \T^n \times J_{\upsilon'} \to \R  \hspace{1mm}| \hspace{1mm} a \in {\bar S}^{\upsilon'}_{\sigma, 2} \hspace{1mm} \mbox{and} \hspace{1mm} |a| =|a|^{\upsilon'}_{\sigma + 2, 0} + |\partial_q a|^{\upsilon'}_{\sigma + 1, \lambda} < \infty\Big\}\\
\mathcal{B} &=& \Big\{b : \T^n \times J_{\upsilon'}  \to \R^n  \hspace{1mm}| \hspace{1mm} b \in {\bar S}^{\upsilon'}_{\sigma, 2} , \hspace{1mm} \mbox{and} \hspace{1mm} |b| =|b|^{\upsilon'}_{\sigma + 2, \lambda} < \infty\Big\}\\
\mathcal{U} &=& \Big\{u :\T^n \times J_{\upsilon'}  \to \R^n  \hspace{1mm}| \hspace{1mm} u,\left( \nabla u\right) \bar W \in {S}^{\upsilon'}_\sigma \\
&& \mbox{and} \hspace{1mm} |u| = \max\{|u|^{\upsilon'}_{\sigma, \lambda}, | \left(\nabla u \right) \bar W|^{\upsilon'}_{\sigma, \lambda} \} < \infty\Big\}\\
\mathcal{V} &=& \Big\{v :\T^n \times J_{\upsilon'}  \to \R^n  \hspace{1mm}| \hspace{1mm} v,\left( \nabla v\right) \bar W \in {S}^{\upsilon'}_\sigma\\
&& \mbox{and} \hspace{1mm} |v| = \max\{|v|^{\upsilon'}_{\sigma, \lambda}, |\left( \nabla v \right) \bar W|^{\upsilon'}_{\sigma, \lambda} \} < \infty\Big\}\\
\mathcal{Z} &=& \Big\{z : \T^n \times J_{\upsilon'}  \to \R^n  \hspace{1mm}| \hspace{1mm} z \in {S}^{\upsilon'}_\sigma , \hspace{1mm} \mbox{and} \hspace{1mm} |z| =|z|^{\upsilon'}_{\sigma, \lambda} < \infty\Big\}\\\mathcal{G} &=& \Big\{g : \T^n \times J_{\upsilon'} \to \R  \hspace{1mm}| \hspace{1mm} g \in {S}^{\upsilon'}_\sigma \hspace{1mm} \mbox{and} \hspace{1mm} |g| =|g|^{\upsilon'}_{\sigma, \lambda} < \infty\Big\}\\
\end{eqnarray*} 
where we recall that the norm $|\cdot |^{\upsilon'}_{\sigma, \lambda}$ is defined by~\eqref{normNC}, while the spaces $S^{\upsilon'}_\sigma$ and ${\bar S}^{\upsilon'}_{\sigma, 2}$ are introduced in Definition \ref{S} and Definition \ref{S2}, respectively. Let $M_n$ be the set of the $n$-dimensional matrices. We introduce two other Banach spaces $\left(\mathcal{M}, |\cdot |\right)$ and $\left(\mathcal{W}, |\cdot |\right)$
\vspace{5mm}
\begin{eqnarray*}
\mathcal{M} &=& \Big\{m : \T^n \times B \times J_{\upsilon'} \to M_n  \hspace{1mm}| \hspace{1mm} m \in \mathcal{\bar S}^{\upsilon'}_{\sigma, 2} \hspace{1mm} \mbox{and} \hspace{1mm} |m| =|m|^{\upsilon'}_{\sigma+2, 1} \le \Upsilon\Big\}\\
\mathcal{W} &=& \Big\{ W:\T^n \to \R^n \hspace{1mm}| \hspace{1mm} W\in C^{\sigma+2}(\T^n)\hspace{1mm} \mbox{and} \hspace{1mm} |W| = |W|_{C^{\sigma+2}}<\infty \Big\}.
\end{eqnarray*} 
where $\Upsilon$ is the positive parameter in~\eqref{H1NC}. We proved in~\cite{Sca22a} that the previous spaces are Banach spaces.

Now, we can correctly define the previous functional $\mathcal{F}$.  Let $\mathcal{F}$ be the following functional
\begin{equation*}
\mathcal{F} :  \mathcal{A} \times \mathcal{B} \times \mathcal{M} \times \mathcal{M} \times \mathcal{W} \times  \mathcal{U} \times  \mathcal{V} \longrightarrow \mathcal{Z} \times \mathcal{G}\\
\end{equation*}
\begin{equation*}
\mathcal{F}(a, b,m , \bar m, W, u, v) = (F_1(b, \bar m, W, u,v), F_2(a,b,m,W,u,v))
\end{equation*}
such that
\begin{eqnarray*}
F_1(b, \bar m,W,u,v) &=& W \circ  (\mathrm{id} + u) - W + b \circ \tilde u + \left(\bar m \circ \tilde \varphi  \right) v - \left(\nabla u \right)\bar W \\,
F_2(a,b, m,W,u,v) &=& \partial_q a\circ \tilde u  + \left(\partial_q W \circ (\mathrm{id} + u) + \partial_q b \circ \tilde u \right) v\\ 
&+& \left(\partial_q m \circ \tilde  \varphi   \right) \cdot v^2 + \left(\nabla v \right)\bar W.
\end{eqnarray*}

\subsection{Homological equation}\label{SecHENC}
We recall some fundamental well-known Gronwall-type inequalities.
\begin{proposition}
\label{Gronwall}
Let $J$ be an interval in $\R$, $t_0 \in J$, and $a$, $b$, $u \in  C(J)$ continuous positive functions. If we assume that 
\begin{equation*}
u(t) \le a(t) + \left|\int_{t_0}^t b(s) u(s) ds \right|, \quad \forall t \in J
\end{equation*}
then it follows that 
\begin{equation}
\label{G1NC}
u(t) \le a(t) + \left|\int_{t_0}^t a(s)b(s) e^{\left|\int_s^t b(\tau)d\tau \right|} ds \right|, \quad \forall t \in J.
\end{equation}
If $a$ is a monotone increasing function and we assume that 
\begin{equation*}
u(t) \le a(t) + \int_{t_0}^t b(s) u(s) ds  \quad \forall t \ge t_0,
\end{equation*}
then we obtain the estimate
\begin{equation}
\label{G2NC}
u(t) \le a(t) e^{\int_{t_0}^t b(s) ds}, \quad \forall t\ge t_0.
\end{equation}
\end{proposition}
\begin{proof}
We refer to~\cite{Amann} for the proof.
\end{proof}

Given $\sigma \ge 1$, $\lambda > 0$ and $\upsilon \ge 0$, this section aims to solve the following equation for the unknown $\varkappa :\T^n \times J_\upsilon \to \R^n$ 
\begin{equation}
\label{HENC}
\begin{cases}
 \partial_q\varkappa(q,t)W(q) + \partial_t \varkappa(q,t) \pm \partial_q W(q) \varkappa(q,t) = z(q,t)\\
W\in C^{\sigma +1}(\T^n), \quad z \in \mathcal{S}^\upsilon_\sigma,\\
|z|^\upsilon_{\sigma, \lambda}<\infty.
\end{cases}
\tag{$HE$}
\end{equation}
If $W(q) \equiv W \in \R^n$ is constant, then the latter translates into the following easier problem
\begin{equation*}
\begin{cases}
 \partial_q\varkappa(q,t)W + \partial_t \varkappa(q,t)  = z(q,t)\\
 z \in \mathcal{S}^\upsilon_\sigma, \quad |z|^\upsilon_{\sigma, \lambda}<\infty.
\end{cases}
\end{equation*}
This is a particular case of the homological equation solved in~\cite{Sca22a}. In this case, we showed that we do not need exponential decay in time. 

Concerning the equation~\eqref{HENC}, we begin by proving several estimates. In what follows, we will widely use the properties of the Hölder norms. For this reason, let us recall the following propositions. 

\begin{proposition}
\label{convexity}
For all $f\in C^{\sigma_1}(\R^n)$, then 
\begin{equation*}
|f|^{\sigma_1-\sigma_0}_{C^\sigma}\le C(\sigma_1) |f|^{\sigma_1-\sigma}_{C^{\sigma_0}}|f|^{\sigma-\sigma_0}_{C^{\sigma_1}} \hspace{3mm} \mbox{for all $0 \le \sigma_0 \le \sigma \le \sigma_1$}.
\end{equation*}
\end{proposition}
\begin{proof}
We refer to~\cite{Hor76} for the proof. 
\end{proof}

Let $D$ be an open subset of $\R^n$, we have the following proposition.
\begin{proposition}
\label{Holder2}
We consider $f$, $g \in C^\sigma(D)$ and $\sigma \ge 0$.
\begin{enumerate}
\item For all $\beta \in \N^{n}$, if $|\beta| + s = \sigma$ then  $\left|{\partial^{|\beta|} \over \partial{x_1}^{\beta_1}... \partial{x_n}^{\beta_n}} f \right|_{C^s} \le |f|_{C^\sigma}$.\\
\item  $|fg|_{C^\sigma} \le C(\sigma)\left(|f|_{C^0}|g|_{C^\sigma} + |f|_{C^\sigma}|g|_{C^0}\right)$. 
\end{enumerate}
Now we consider composite functions. Let $z$ be defined on $D_1 \subset \R^n$ and takes its values on $D_2 \subset \R^n$ where $f$ is defined. 

If $\sigma < 1$, $f \in C^1(D_2)$, $z \in C^\sigma (D_1)$ then $f\circ z \in C^\sigma(D_1)$ 
\begin{enumerate}
\item[3.] $|f \circ z|_{C^\sigma} \le C(|f|_{C^1}|z|_{C^\sigma}+ |f|_{C^0})$.
\end{enumerate}
If $\sigma < 1$, $f \in C^\sigma(D_2)$, $z \in C^1 (D_1)$ then $f\circ z \in C^\sigma(D_1)$  
\begin{enumerate}
\item[4.] $|f \circ z|_{C^\sigma} \le C(|f|_{C^\sigma}|\nabla z|^\sigma_{C^0}+ |f|_{C^0})$.
\end{enumerate}
If $\sigma \ge 1$ and $f \in C^\sigma (D_2)$, $z \in C^\sigma (D_1)$ then $f\circ z \in C^\sigma(D_1)$ 
\begin{enumerate}
\item[5.] $|f \circ z|_{C^\sigma} \le C(\sigma) \left(|f|_{C^\sigma}|\nabla z|^\sigma_{C^0} + |f|_{C^1}|\nabla z|_{C^{\sigma-1}}+ |f|_{C^0}\right)$.
\end{enumerate}
\end{proposition}
\begin{proof}
We refer to~\cite{Hor76} for the proof of \textit{1.} and \textit{2.}. The following two (\textit{3.} and \textit{4.}) are obvious, and we refer to~\cite{Sca22a} for the last one.
\end{proof}

Let $\phi_W^t$ be the flow at time $t$ of $W(q)$. As usual, $C(\cdot)$ stands for constants depending on $n$ and the other parameters into brackets. On the other hand, $C$ means constants depending on $n$.
\begin{lemma}
\label{psiNC}
For all $t \in \R$ 
\begin{equation}
\label{estimatepsiposNC}
|\partial_q \phi_W^t|_{C^{\sigma-1}} \le C(\sigma)\left(1 + |\partial_q W|_{C^{\sigma-1}}  |t|\right) e^{c_\sigma |\partial_q W|_{C^0}|t|},
\end{equation}
with a positive constant $c_\sigma \ge 1$ depending on $n$ and $\sigma$. 
\end{lemma}
By~\eqref{estimatepsiposNC}, we note that when $\sigma = 1$ and $t \in \R$
\begin{eqnarray*}
|\partial_q \phi_W^t|_{C^0} &\le& C\left(1 + |\partial_q W|_{C^0}  |t|\right) e^{c_1 |\partial_q W|_{C^0}|t|}\le C e^{\bar c_1 |\partial_q W|_{C^0}|t|}
\end{eqnarray*}
for a suitable $\bar c_1 >c_1$. 
\begin{proof}
For all $q \in \T^n$, by the fundamental theorem of calculus, we can write $\phi_W^t$ in the following form
\begin{equation*}
\phi_W^t(q) = q + \int_0^t W \circ \phi_W^\tau(q)d\tau.
\end{equation*}
Therefore, taking the derivative with respect to $q$ 
\begin{equation*}
\partial_q \phi_W^t(q) = \mathrm{Id} + \int_0^t \partial_qW \circ \phi_W^\tau(q)\partial_q\phi_W^\tau(q)d\tau,
\end{equation*}
where $\mathrm{Id}$ stands for the identity matrix. 
We assume $t \ge 0$. Then, we can estimate the norm $C^{\sigma-1}$ of the left-hand side of the latter as follows
\begin{equation}
\label{entermediatestimatephiNC}
|\partial_q \phi_W^t|_{C^{\sigma-1}} \le 1 + \int_0^t|\partial_qW \circ \phi_W^\tau\partial_q\phi_W^\tau|_{C^{\sigma-1}} d\tau. 
\end{equation}

\vspace{5mm}
\textit{Case $\sigma=1$}. By Proposition \ref{Holder2}
\begin{equation*}
|\partial_q \phi_W^t|_{C^0} \le 1 + C\int_0^t |\partial_q W|_{C^0}|\partial_q \phi^\tau_W|_{C^0}d\tau,
\end{equation*}
for a suitable constant $C$. Then, thanks to~\eqref{G2NC}
\begin{equation}
\label{phi1NC}
|\partial_q \phi_W^t|_{C^0} \le e^{c_1|\partial_q W|_{C^0}t}
\end{equation}
for a suitable constant $c_1 \ge 1$.

It remains to verify~\eqref{estimatepsiposNC} when $\sigma >1$. By Proposition~\eqref{Holder2}, we can estimate the norm on the right-hand side of~\eqref{entermediatestimatephiNC} as follows
\begin{eqnarray*}
|\partial_qW \circ \phi_W^\tau\partial_q\phi_W^\tau|_{C^{\sigma-1}} &\le& C(\sigma)\left(|\partial_qW \circ \phi_W^\tau|_{C^{\sigma-1}} | \partial_q\phi_W^\tau|_{C^0} +  |\partial_qW|_{C^0}|\partial_q\phi_W^\tau|_{C^{\sigma-1}}\right).
\end{eqnarray*}
Hence, we can rewrite~\eqref{estimatepsiposNC} in the following form
\begin{eqnarray}
|\partial_q \phi_W^t|_{C^{\sigma-1}} &\le& 1 + C(\sigma)\int_0^t |\partial_qW \circ \phi_W^\tau|_{C^{\sigma-1}} | \partial_q\phi_W^\tau|_{C^0}d\tau \nonumber\\
\label{entermediatestimatephi2NC}
&+&C(\sigma)\int_0^t|\partial_qW|_{C^0}|\partial_q\phi_W^\tau|_{C^{\sigma-1}}, d\tau. 
\end{eqnarray} 
We need to treat cases $1<\sigma <2$ and $\sigma \ge 2$ separately.

\vspace{5mm}
\textit{Case $1<\sigma <2$}. Thanks to Proposition~\eqref{Holder2}, in particular property \textit{4.},  
\begin{eqnarray*}
 |\partial_qW \circ \phi_W^\tau|_{C^{\sigma-1}} \le C(\sigma)\left(|\partial_q W|_{C^{\sigma-1}}|\partial_q \phi^\tau_W|^{\sigma-1}_{C^0} + |\partial_q W|_{C^0}\right).
\end{eqnarray*}
Replacing the latter into~\eqref{entermediatestimatephi2NC}, we can rewrite it in the following way
\begin{eqnarray*}
|\partial_q \phi_W^t|_{C^{\sigma-1}} &\le& 1 + C(\sigma)\int_0^t |\partial_q W|_{C^0}|\partial_q \phi^\tau_W|_{C^0}d\tau\\
&+& C(\sigma)\int_0^t |\partial_q W|_{C^{\sigma-1}}|\partial_q \phi^\tau_W|^\sigma_{C^0}d\tau + C(\sigma)\int_0^t|\partial_qW|_{C^0}|\partial_q\phi_W^\tau|_{C^{\sigma-1}} d\tau. 
\end{eqnarray*} 
Now,~\eqref{phi1NC} allows us to find an upper bound for the first two integrals on the right-hand side of the latter
\begin{eqnarray*}
\int_0^t |\partial_q W|_{C^0}|\partial_q \phi^\tau_W|_{C^0}d\tau &\le& |\partial_q W|_{C^0}\int_0^t e^{c_1|\partial_q W|_{C^0}\tau}d\tau = {e^{c_1 |\partial_q W|_{C^0}t}-1 \over c_1}\\
\int_0^t |\partial_q W|_{C^{\sigma-1}}|\partial_q \phi^\tau_W|^\sigma_{C^0}d\tau &\le& |\partial_q W|_{C^{\sigma-1}}\int_0^t e^{c_1\sigma|\partial_q W|_{C^0}\tau}d\tau \le |\partial_q W|_{C^{\sigma-1}} t e^{c_1\sigma|\partial_q W|_{C^0}t}.
\end{eqnarray*} 
In the second line of the latter, rather than calculating the integral, we prefer using the trivial estimate $e^{c_1\sigma|\partial_q W|_{C^0}\tau} \le e^{c_1\sigma|\partial_q W|_{C^0}t}$ to avoid a division by $|\partial_q W|_{C^0}$ since we do not assume it is not zero.

Hence, we can estimate $|\partial_q \phi_W^t|_{C^{\sigma-1}} $ as follows
\begin{eqnarray*}
|\partial_q \phi_W^t|_{C^{\sigma-1}} &\le& 1 + C(\sigma){e^{c_1 |\partial_q W|_{C^0}t}-1 \over c_1} + C(\sigma) |\partial_q W|_{C^{\sigma-1}} t e^{c_1\sigma|\partial_q W|_{C^0}t}\\
&+& C(\sigma) \int_0^t|\partial_qW|_{C^0}|\partial_q\phi_W^\tau|_{C^{\sigma-1}}, d\tau,\\
&\le& C(\sigma)\left(1  +  |\partial_q W|_{C^{\sigma-1}} t \right)e^{c_1\sigma|\partial_q W|_{C^0}t} + C(\sigma) \int_0^t|\partial_qW|_{C^0}|\partial_q\phi_W^\tau|_{C^{\sigma-1}}, d\tau.
\end{eqnarray*} 
Then, thanks to the Gronwall inequality~\eqref{G2NC}
\begin{eqnarray*}
|\partial_q \phi_W^t|_{C^{\sigma-1}} &\le&C(\sigma)\left(1  +  |\partial_q W|_{C^{\sigma-1}} t \right)e^{c_1\sigma|\partial_q W|_{C^0}t}e^{C(\sigma) \int_0^t |\partial_q W|_{C^0}d\tau}\\
&\le& C(\sigma)\left(1  +  |\partial_q W|_{C^{\sigma-1}} t \right)e^{ c_\sigma|\partial_q W|_{C^0}t}
\end{eqnarray*} 
for a suitable constant $c_\sigma \ge c_1 \sigma$. This concludes the proof for the case $1<\sigma <2$. The general case $\sigma >2$ is quite similar to the previous one. The main difference lies on the estimation of $ |\partial_qW \circ \phi_W^\tau|_{C^{\sigma-1}} $. 

\vspace{5mm}
\textit{Case $\sigma >2$}. By Proposition \ref{Holder2}, especially property \textit{5.},
\begin{eqnarray*}
 |\partial_qW \circ \phi_W^\tau|_{C^{\sigma-1}} &\le& C(\sigma)\Big(|\partial_q W|_{C^{\sigma-1}}|\partial_q \phi^\tau_W|^{\sigma-1}_{C^0} + |\partial_q W|_{C^1}|\partial_q \phi^\tau_W|_{C^{\sigma-2}}  + |\partial_q W|_{C^0}\Big)
\end{eqnarray*}
and replacing the latter into~\eqref{entermediatestimatephi2NC}, we can estimate $|\partial_q \phi_W^t|_{C^{\sigma-1}}$ as follows
\begin{eqnarray*}
|\partial_q \phi_W^t|_{C^{\sigma-1}} &\le& 1 + C(\sigma)\int_0^t |\partial_q W|_{C^0}|\partial_q \phi^\tau_W|_{C^0}d\tau + C(\sigma)\int_0^t |\partial_q W|_{C^{\sigma-1}}|\partial_q \phi^\tau_W|^\sigma_{C^0}d\tau\\
&+& C(\sigma)\int_0^t |\partial_q W|_{C^1}|\partial_q \phi^\tau_W|_{C^{\sigma-2}}|\partial_q \phi^\tau_W|_{C^0}d\tau\\
&+& C(\sigma)\int_0^t|\partial_qW|_{C^0}|\partial_q\phi_W^\tau|_{C^{\sigma-1}} d\tau. 
\end{eqnarray*} 
We have already estimated the first two integrals on the right-hand side of the latter. It remains the integral in the second line. By the convexity property of the Hölder norms (Proposition \ref{convexity}), for all fixed $\tau$
\begin{eqnarray*}
|\partial_q W|_{C^1}|\partial_q \phi^\tau_W|_{C^{\sigma-2}} &\le& C(\sigma) \left(|\partial_q W|^{\sigma-2 \over \sigma-1}_{C^0}|\partial_q W|_{C^{\sigma-1}}^{1\over \sigma -1}\right) \left(|\partial_q \phi^\tau_W|^{1\over \sigma -1}_{C^0}|\partial_q \phi^\tau_W|_{C^{\sigma-1}}^{\sigma-2\over \sigma -1}\right) 
\end{eqnarray*}
and hence
\begin{eqnarray*}
|\partial_q W|_{C^1}|\partial_q \phi^\tau_W|_{C^{\sigma-2}}|\partial_q \phi^\tau_W|_{C^0} &\le&C(\sigma)  \left(|\partial_q W|_{C^0}|\partial_q \phi^\tau_W|_{C^{\sigma-1}}\right)^{\sigma-2 \over \sigma-1}  \left(|\partial_q W|_{C^{\sigma-1}}|\partial_q \phi^\tau_W|_{C^0}^\sigma\right)^{1 \over \sigma-1}.
\end{eqnarray*}
From $a^\lambda b^{1-\lambda} \le C(a+b)$ for $0 < \lambda < 1$, we have that 
\begin{eqnarray*}
|\partial_q W|_{C^1}|\partial_q \phi^\tau_W|_{C^{\sigma-2}}|\partial_q \phi^\tau_W|_{C^0} &\le& C(\sigma) \left(|\partial_q W|_{C^0}|\partial_q \phi^\tau_W|_{C^{\sigma-1}} + |\partial_q W|_{C^{\sigma-1}}|\partial_q \phi^\tau_W|_{C^0}^\sigma\right).
\end{eqnarray*}
Furthermore, replacing the latter in the previous estimate of $|\partial_q \phi_W^t|_{C^{\sigma-1}} $, we obtain 
\begin{eqnarray*}
|\partial_q \phi_W^t|_{C^{\sigma-1}} &\le& 1 + C(\sigma)\int_0^t |\partial_q W|_{C^0}|\partial_q \phi^\tau_W|_{C^0}d\tau\\
&+& C(\sigma)\int_0^t |\partial_q W|_{C^{\sigma-1}}|\partial_q \phi^\tau_W|^\sigma_{C^0}d\tau + C(\sigma)\int_0^t|\partial_qW|_{C^0}|\partial_q\phi_W^\tau|_{C^{\sigma-1}} d\tau. 
\end{eqnarray*} 
Now, similarly to the previous case ($1<\sigma<2$), we conclude the proof of~\eqref{estimatepsiposNC} also in this general case. 
Similarly, we have the claim when $t \le 0$.
\end{proof}
 Now, we consider $R : \T^n \times J_\upsilon \times J_\upsilon \to M_n$, where $M_n$ is the set of the $n$-dimensional matrices. For all $(q,\tau, t) \in \T^n \times J_\upsilon \times J_\upsilon$, $R(q, t,\tau)$ is the matrix having elements equal to $r_{ij}(q,t,\tau)$ for all $1 \le i,j \le n$. In other words, $R(q, t,\tau) = \{r_{ij}(q, t,\tau) \}_{1 \le i,j \le n}$. We define the following family of norms 
\begin{equation*}
|R^t_\tau|_{C^s} = \max_{1 \le i,j \le n}|r_{ij}(q,t,\tau)|_{C^s},
\end{equation*}
for positive real parameters $s \ge 0$. We consider the following system that plays an important role in the solution of the homological equation~\eqref{HENC}
\begin{equation}
\label{RNC}
\begin{cases}
\dot R(q, t,\tau) = \mp \partial_q W \circ \phi_W^t(q) R(q, t,\tau) \\
R(q,\tau,\tau) = \mathrm{Id}.
\end{cases}
\tag{R}
\end{equation}
where $W$ is defined in~\eqref{HENC}. For all fixed $\tau$, $t \in J_\upsilon$, in what follows we denote $R_\tau^t(q) = R(q, t, \tau)$. 

\begin{lemma}
\label{LemmaRNC}
The latter system admits a unique solution. Moreover, for all $\tau$, $t \in J_\upsilon$ with $\tau \ge t$, letting $\tilde R(q,t,\tau) = R(\phi_W^{-\tau}(q),t,\tau)$, we have the following estimates 
\begin{eqnarray}
\label{stimeR0NC}
|R^t_\tau|_{C^0} &\le& e^{c_0^R|\partial_q W|_{C^0}(\tau-t)}\\
\label{stimeRNC}
|\tilde R^t_\tau|_{C^\sigma} &\le& C(\sigma)\left(1 + |\partial_q W|_{C^\sigma} (\tau - t)\right)e^{c^R_\sigma |\partial_q W|_{C^0} (\tau-t)}\\
&+&C(\sigma)|\partial_q W|_{C^1}|\partial_q W|_{C^{\sigma-1}} (\tau - t)^2 e^{c^R_\sigma |\partial_q W|_{C^0} (\tau-t)}\nonumber
\end{eqnarray}
with positive constants $c_0^R>0$ and  $c^R_\sigma \ge c_\sigma$. We point out that $c_\sigma$ is the positive constant introduced in the previous lemma. 
\end{lemma}
Before the proof, we observe that when $\sigma = 1$, thanks to~\eqref{stimeRNC}, 
\begin{eqnarray*}
|\tilde R^t_\tau|_{C^1} &\le& C\left(1 + |\partial_q W|_{C^1} (\tau - t) + |\partial_q W|_{C^1}|\partial_q W|_{C^0} (\tau - t)^2 \right)e^{c^R_1 |\partial_q W|_{C^0} (\tau-t)}\\
&\le& C\left(1 + |\partial_q W|_{C^1} (\tau - t)\right)e^{c^R_1 |\partial_q W|_{C^0} (\tau-t)}\\
&+&C|\partial_q W|_{C^1}(\tau - t) e^{|\partial_q W|_{C^0} (\tau-t)}e^{c^R_1 |\partial_q W|_{C^0} (\tau-t)}\\
&\le& C\left(1 + |\partial_q W|_{C^1} (\tau - t)\right)e^{\bar c^R_1 |\partial_q W|_{C^0} (\tau-t)}\\
\end{eqnarray*}
for a suitable $\bar c_1^R > c_1^R$.
\begin{proof}
We prove this lemma in the case $\dot R(q, t,\tau) = \partial_q W \circ \phi_W^t(q) R(q, t,\tau)$. The other case ($\dot R(q, t,\tau) =- \partial_q W \circ \phi_W^t(q) R(q, t,\tau)$) can be proved similarly.
For all $q \in \T^n$, a unique solution of~\eqref{RNC} exists by the existence and uniqueness theorem. It remains to prove the estimates. 

By the fundamental theorem of calculus, we can write $R$ as follows
\begin{equation}
\label{formaRNC}
R^t_\tau(q) = \mathrm{Id} - \int_t^\tau \partial_q W \circ \phi^s_W(q) R^s_\tau(q) ds.
\end{equation}
for all $q \in \T^n$ and $t$, $\tau \in J_\upsilon$ with $\tau \ge t$. 
Then, thanks to the latter, we can estimate $|R^t_\tau|_{C^0}$ in the following way
\begin{equation*}
|R^t_\tau|_{C^0} \le 1 + C\int_t^\tau |\partial_qW|_{C^0} |R^s_\tau|_{C^0}ds.
\end{equation*}
Therefore, by the Gronwall inequality~\eqref{G2NC}, we have
\begin{equation*}
|R^t_\tau|_{C^0} \le e^{\int_t^\tau C|\partial_qW|_{C^0}ds} \le e^{c_R^0 |\partial_qW|_{C^0}(\tau-t)}
\end{equation*}
for a suitable positive constant $c_0^R$. This concludes the proof of~\eqref{stimeR0NC}. Now, we prove~\eqref{stimeRNC}. By~\eqref{formaRNC}, we can write $\tilde R^t_\tau$ in the following form
\begin{equation*}
\tilde R^t_\tau(q) = R^t_\tau \circ \phi_W^{-\tau}(q) = \mathrm{Id} - \int_t^\tau \partial_q W \circ \phi^{s-\tau}_W(q) \tilde R^s_\tau(q) ds.
\end{equation*}
Hence, we can estimate $|\tilde R^t_\tau|_{C^\sigma}$ in such a way that
\begin{equation}
\label{partialestR}
|\tilde R^t_\tau|_{C^\sigma}\le 1 + \int_t^\tau |\partial_q W \circ \phi^{s-\tau}_W \tilde R^s_\tau|_{C^\sigma}ds.
\end{equation}
We will estimate the norm into the integral on the right-hand side of the latter using Proposition \ref{Holder2} and~\eqref{stimeR0NC}. The claim is a consequence of the Gronwall inequality~\eqref{G1NC}. 

As a consequence of Proposition \ref{Holder},
\begin{eqnarray*}
|\partial_q W  \circ \phi_W^{s-\tau} \tilde R^s_\tau|_{C^\sigma} &\le& C(\sigma) \left(|\partial_q W  \circ \phi_W^{s-\tau}|_{C^\sigma}|R^s_\tau|_{C^0} + |\partial_qW|_{C^0}|\tilde R^s_\tau|_{C^\sigma}\right)\\
|\partial_q W  \circ \phi_W^{s-\tau}|_{C^\sigma} &\le& C(\sigma) \left(|\partial_q W|_{C^\sigma}|\partial_q \phi_W^{s-\tau}|^\sigma_{C^0} + |\partial_q W|_{C^1}|\partial_q \phi_W^{s-\tau}|_{C^{\sigma-1}} + |\partial_qW|_{C^0} \right)
\end{eqnarray*}
and replacing the latter into~\eqref{partialestR}
\begin{eqnarray*}
|\tilde R^t_\tau|_{C^\sigma} &\le& 1 + C(\sigma)\int_t^\tau|\partial_q W|_{C^\sigma}|\partial_q \phi_W^{s-\tau}|^\sigma_{C^0}  |R^s_\tau|_{C^0}ds\\
&+& C(\sigma)\int_t^\tau |\partial_q W|_{C^1}|\partial_q \phi_W^{s-\tau}|_{C^{\sigma-1}} |R^s_\tau|_{C^0}ds\\
&+& C(\sigma)\int_t^\tau |\partial_q W|_{C^0} |R^s_\tau|_{C^0}ds + C(\sigma)\int_t^\tau |\partial_q W|_{C^0} |\tilde R^s_\tau|_{C^\sigma}ds.
\end{eqnarray*}
Now, by~\eqref{stimeR0NC} and Lemma \ref{psiNC}, we can estimate the first three integrals on the right-hand side of the latter
\begin{eqnarray*}
\int_t^\tau|\partial_q W|_{C^\sigma}|\partial_q \phi_W^{s-\tau}|^\sigma_{C^0}  |R^s_\tau|_{C^0}ds &\le& C(\sigma)|\partial_q W|_{C^\sigma}\int_t^\tau e^{\bar c_1 \sigma|\partial_q W|_{C^0}(\tau -s)} e^{c_0^R|\partial_q W|_{C^0}(\tau - s)}ds \\
&\le& C(\sigma)|\partial_q W|_{C^\sigma}(\tau - t) e^{(\bar c_1\sigma + c_0^R)|\partial_q W|_{C^0}(\tau -t)}\\
\int_t^\tau |\partial_q W|_{C^1}|\partial_q \phi_W^{s-\tau}|_{C^{\sigma-1}} |R^s_\tau|_{C^0}ds &\le&C(\sigma)  |\partial_q W|_{C^1}\int_t^\tau e^{(c_\sigma + c_0^R) |\partial_q W|_{C^0}(\tau-s)} ds \\
&+&C(\sigma)  |\partial_q W|_{C^1} |\partial_q W|_{C^{\sigma -1}}\int_t^\tau (\tau - s)e^{(c_\sigma + c_0^R) |\partial_q W|_{C^0}(\tau-s)} ds \\
&\le&C(\sigma) |\partial_q W|_{C^1}(\tau - t) e^{(c_\sigma + c_0^R) |\partial_q W|_{C^0}(\tau-t)} \\
&+&C(\sigma)  |\partial_q W|_{C^1} |\partial_q W|_{C^{\sigma -1}}(\tau - t)^2 e^{(c_\sigma + c_0^R) |\partial_q W|_{C^0}(\tau-t)} \\
\int_t^\tau |\partial_q W|_{C^0} |R^s_\tau|_{C^0}ds &\le&  |\partial_q W|_{C^0}\int_t^\tau e^{c_0^R |\partial_q W|_{C^0}(\tau - s)}ds\\
&=& {1 \over c_0^R} \left(e^{c_0^R |\partial_q W|_{C^0}(\tau-t)} - 1\right).
\end{eqnarray*}
Similarly to the previous lemma, in the first two integrals on the left-hand side of the latter, we use some trivial inequalities to avoid the division by $ |\partial_q W|_{C^0}$.

Now, we recall that $c_\sigma \ge c_1 \sigma$. Moreover, we can assume without loss of generality that $c_\sigma \ge \max\{c_1 \sigma, \bar c_1 \sigma\}$.
Hence, by the above estimations 
\begin{eqnarray*}
|\tilde R^t_\tau|_{C^\sigma} &\le& 1 +  C(\sigma)\left(e^{c_0^R |\partial_q W|_{C^0}(\tau-t)} - 1\right) + C(\sigma)|\partial_q W|_{C^\sigma}(\tau - t) e^{(\bar c_1\sigma + c_0^R)|\partial_q W|_{C^0}(\tau -t)}\\
&+& C(\sigma) |\partial_q W|_{C^1}(\tau - t) e^{(c_\sigma + c_0^R) |\partial_q W|_{C^0}(\tau-t)} \\
&+&C(\sigma)  |\partial_q W|_{C^1} |\partial_q W|_{C^{\sigma -1}}(\tau - t)^2 e^{(c_\sigma + c_0^R) |\partial_q W|_{C^0}(\tau-t)}\\
&+&  C(\sigma)\int_t^\tau |\partial_q W|_{C^0} |\tilde R^s_\tau|_{C^\sigma}ds\\
&\le& C(\sigma) \left(1 + |\partial_q W|_{C^\sigma}(\tau - t)\right)e^{(c_\sigma + c_0^R) |\partial_q W|_{C^0}(\tau-t)} \\
&+&C(\sigma)  |\partial_q W|_{C^1} |\partial_q W|_{C^{\sigma -1}}(\tau - t)^2 e^{(c_\sigma + c_0^R) |\partial_q W|_{C^0}(\tau-t)}\\
&+& \left|C(\sigma)\int_\tau^t |\partial_q W|_{C^0} |\tilde R^s_\tau|_{C^\sigma}ds\right|.
\end{eqnarray*}
At this moment, we define the following function 
\begin{eqnarray*}
a(t) &=& C(\sigma) \left(1 + |\partial_q W|_{C^\sigma}(\tau - t)\right)e^{(c_\sigma + c_0^R) |\partial_q W|_{C^0}(\tau-t)}\\
&+& C(\sigma)  |\partial_q W|_{C^1} |\partial_q W|_{C^{\sigma -1}}(\tau - t)^2 e^{(c_\sigma + c_0^R) |\partial_q W|_{C^0}(\tau-t)}
\end{eqnarray*}
and we rewrite the latter in the following way
\begin{equation*}
|\tilde R^t_\tau|_{C^\sigma} \le a(t) + \left|C(\sigma)\int_\tau^t |\partial_q W|_{C^0} |\tilde R^s_\tau|_{C^\sigma}ds\right|.
\end{equation*}

However, it is straightforward to verify that $a$ is a monotone decreasing function. Hence, by the more general inequality~\eqref{G1NC}
\begin{eqnarray*}
|\tilde R^t_\tau|_{C^\sigma} &\le& a(t) +C(\sigma)\left|\int_\tau^t a(s)  |\partial_q W|_{C^0} e^{\left|C(\sigma)\int^t_s  |\partial_q W|_{C^0}d\delta\right|}ds\right|\\
&\le& a(t)\left(1 +C(\sigma) |\partial_q W|_{C^0}  \int^\tau_t  e^{C(\sigma) |\partial_q W|_{C^0}(s-t)}ds \right)\\
&=& a(t)\left(1 + \left(e^{C(\sigma) |\partial_q W|_{C^0}(\tau-t)}-1\right) \right)\\
&\le& a(t)\left(1 + e^{C(\sigma) |\partial_q W|_{C^0}(\tau-t)}\right)\\
&\le& C(\sigma) \left(1 + |\partial_q W|_{C^\sigma}(\tau - t)\right)e^{c_\sigma^R |\partial_q W|_{C^0}(\tau-t)}\\
&+& C(\sigma)  |\partial_q W|_{C^1} |\partial_q W|_{C^{\sigma -1}}(\tau - t)^2 e^{c_\sigma^R |\partial_q W|_{C^0}(\tau-t)}
\end{eqnarray*}
for a suitable constant $c_\sigma^R \ge c_\sigma + c_0^R$. 
\end{proof}
We observe that the constant $c_\sigma^R$, as for $c_\sigma$, goes to infinity if $\sigma \to \infty$. This means that, in order to solve the homological equation, we must counter the growth of $c_\sigma^R$ and $c_\sigma$ assuming $\lambda$ sufficiently large. 

\begin{lemma}[\textbf{Homological equation}]
\label{homoeqlemmaNC} 
There exists a solution $\varkappa$, $\left(\nabla \varkappa\right)\bar W \in \mathcal{S}^\upsilon_\sigma$ of~\eqref{HENC}. Moreover, letting $c_\sigma^\varkappa = \max\{c_\sigma + c_0^R, c_\sigma^R + c_0^R, \bar c_1^R + c_\sigma\}$, if
\begin{equation}
\label{HElambdaNC}
\lambda > c_\sigma^\varkappa |\partial_q W|_{C^0} 
\end{equation}
then, 
\begin{eqnarray}
\label{varkappaNC}
|\varkappa|^\upsilon_{\sigma,\lambda} &\le& C(\sigma) \left({1 \over \lambda - c_\sigma^\varkappa |\partial_q W|_{C^0}}\right)|z|^\upsilon_{\sigma, \lambda}\\
&+&C(\sigma)\left({|\partial_q W|_{C^\sigma} \over \left(\lambda - c_\sigma^\varkappa |\partial_q W|_{C^0}\right)^2} + {|\partial_q W|_{C^1}  |\partial_q W|_{C^{\sigma-1}} \over \left(\lambda - c_\sigma^\varkappa |\partial_q W|_{C^0}\right)^3}\right) |z|^\upsilon_{\sigma, \lambda}. \nonumber
\end{eqnarray}
\end{lemma}
\begin{proof}
\textit{Existence}: Let us define the following transformation 
\begin{align}
\label{hNC}
&h:\T^n \times J_\upsilon \longrightarrow \T^n \times J_\upsilon\\
&h(q,t) = (\phi_W^{-t}(q) , t)\nonumber
\end{align}
where $\phi^t_W$ is the flow of $W$ previously introduced. We claim that it suffices to prove the first part of this lemma for the much easier equation
\begin{equation}
\label{HE2NC}
\partial_t \kappa(q,t) \pm \partial_q W \circ \phi_W^t(q) \kappa(q,t) = z \circ h^{-1}(q,t).
\end{equation} 
If $\kappa$ is a solution of the latter, then $\varkappa = \kappa \circ h$ is a solution of~\eqref{HENC} and viceversa. We prove this claim. Let $\varkappa$ be a solution of~\eqref{HENC}, 
\begin{eqnarray*}
\partial_t \left(\varkappa \circ h^{-1}\right) \pm \left(\partial_q W\circ \phi^t_W \right)\left(\varkappa \circ h^{-1}\right) &=& \left(\partial_q\varkappa \circ h^{-1}\right)\dot \phi^t_W + \partial_t \varkappa \circ h^{-1} \\
&\pm& \left(\partial_q W\circ \phi^t_W \right)\left(\varkappa \circ h^{-1}\right)\\
&=& \left(\partial_q\varkappa \circ h^{-1}\right)\left(W \circ \phi^t_W\right) + \partial_t \varkappa \circ h^{-1}\\
&\pm& \left(\partial_q W\circ \phi^t_W \right)\left(\varkappa \circ h^{-1}\right)\\
&=& z \circ h^{-1}
\end{eqnarray*}
where $\dot \phi^t_W$ stands for the derivative of $\phi^t_W$ with respect to $t$. It is obviously equal to $W \circ \phi^t_W$. Furthermore, the last equality is a consequence of~\eqref{HENC}. This proves that  $\varkappa \circ h^{-1}$ is a solution of~\eqref{HE2NC}. Let us first show 
\begin{equation}
\label{HEprimaparte1}
\partial_q \phi^{-t}_W(q) W(q) = \dot \phi^{-t}_W(q)
\end{equation} 
for all $(q,t) \in \T^n \times J_\upsilon$. We know that $\phi^t_W$ is the flow of $W$, and hence  $\dot \phi^{-t}_W = W \circ \phi^{-t}_W $. This implies that it is enough to prove the equivalent equality
\begin{equation}
\label{HEprimaparte2}
W \circ \phi^{-t}_W(q) = \partial_q \phi^{-t}_W(q) W(q)
\end{equation}
for all $(q,t) \in \T^n \times J_\upsilon$.
We observe that the pull-back of $W$ by $ \phi^{-t}_W$ is equal to $W$. In others words, $\left(\phi^{-t}_W\right)^*W = W$ where $W = \left(\phi^{-t}_W\right)^*W = \left(\partial_q \phi^{-t}_W\right)^{-1}W \circ \phi^{-t}_W$. This proves~\eqref{HEprimaparte2} and hence~\eqref{HEprimaparte1}.

Now, let $\kappa$ be a solution of~\eqref{HE2NC}, then 
\begin{eqnarray*}
\partial_q \left(\kappa \circ h\right)W + \partial_t\left(\kappa \circ h\right) \pm \partial_q W\left(\kappa \circ h\right) &=& \left(\partial_q\kappa \circ h\right)\partial_q \phi^{-t}_W W- \left(\partial_q\kappa \circ h\right)\dot \phi^{-t}_W\\
&+& \partial_t\kappa \circ h \pm \partial_q W\left(\kappa \circ h\right)\\
&=& \partial_t\kappa \circ h \pm \partial_q W\left(\kappa \circ h\right) = z.
\end{eqnarray*}
Hence $\kappa \circ h$ is a solution of~\eqref{HENC}, where the last line of the latter is a consequence of~\eqref{HEprimaparte1} and the last equality of~\eqref{HE2NC}. This proves the claim. 

For all $q \in \T^n$, let $R(q, t ,\upsilon)$ be the unique solution of~\eqref{RNC}. For all $(q, t) \in \T^n \times J_\upsilon$ a solution $\kappa$ of~\eqref{HE2NC} exists and 
\begin{eqnarray*}
\kappa (q,t) &=& R(q, t, \upsilon) e(q) - \int_{\upsilon}^t R(q, t, \tau)  z\circ h^{-1}(q, \tau) d\tau\\
&=& R(q, t,\upsilon)\left(e(q) - \int_{\upsilon}^t R(q, \upsilon, \tau) z \circ h^{-1}(q, \tau) d\tau  \right)
\end{eqnarray*}
where $e$ is a function defined on the torus.

\vspace{5mm}
\textit{Estimates}: We choose $e$ equal to 
\begin{equation*}
e(q) = \int_{\upsilon}^{+\infty}R(q, \upsilon, \tau) z \circ h^{-1}(q, \tau) d\tau
\end{equation*}
for all $q \in \T^n$. It is well defined because, by Lemma \ref{LemmaRNC} and~\eqref{HElambdaNC},
\begin{eqnarray*}
\left|\int_{\upsilon}^{+\infty}R(q, \upsilon, \tau) z\circ h^{-1}(q, \tau) d\tau \right| &\le& C\int_{\upsilon}^{+\infty}|R^{\upsilon}_\tau|_{C^0} |z^\tau|_{C^0} d\tau\\
&\le& C |z|^\upsilon_{0, \lambda}\int_{\upsilon}^{+\infty}e^{\left(c_0^R|\partial_q W|_{C^0} -\lambda \right)s}ds\\
&=& {C |z|^\upsilon_{0, \lambda} \over \lambda - c_0^R|\partial_q W|_{C^0}}e^{\left(c_0^R|\partial_q W|_{C^0} -\lambda \right)\upsilon}
\end{eqnarray*}
Therefore,  for all $(q,t) \in \T^n \times J_\upsilon$,
\begin{eqnarray*}
\varkappa(q,t) &=& \kappa \circ h(q,t) = -\int_t^{+\infty} R^t_\tau\circ \phi^{-t}_W(q) z^\tau\circ \phi^{\tau-t}_W(q)d\tau\\
&=&-\int_t^{+\infty} R^t_\tau\circ \phi^{-\tau}_W\circ \phi^{\tau-t}_W(q) z^\tau\circ \phi^{\tau-t}_W(q)d\tau \\
&=& -\int_t^{+\infty} \tilde R^t_\tau\circ \phi^{\tau-t}_W(q) z^\tau\circ \phi^{\tau-t}_W(q)d\tau
\end{eqnarray*}
is the solution of~\eqref{HENC} we are looking for. 

The estimate~\eqref{varkappaNC} is a consequence of Proposition \ref{Holder2}, Lemma \ref{psiNC}, Lemma \ref{LemmaRNC} and~\eqref{HElambdaNC}. For all fixed $t \in J_\upsilon$, by Proposition \ref{Holder2}, we can estimate $|\varkappa^t|_{C^\sigma}$ as follows
\begin{equation*}
|\varkappa^t|_{C^\sigma} \le C(\sigma) \int_t^{+\infty} |\tilde R^t_\tau\circ \phi^{\tau-t}_W|_{C^\sigma}|z^\tau|_{C^0} +  |R^{t}_\tau|_{C^0}|z^\tau\circ \phi^{\tau-t}_W|_{C^\sigma}d\tau.
\end{equation*}
Always using Proposition \ref{Holder2}
\begin{eqnarray*}
|z^\tau\circ \phi^{\tau-t}_W|_{C^\sigma}&\le& C(\sigma)|z^\tau|_{C^\sigma} \left(|\partial_q \phi_W^{\tau-t}|^\sigma_{C^0} + |\partial_q \phi_W^{\tau-t}|_{C^{\sigma-1}} + 1\right)\\
 |\tilde R^t_\tau\circ \phi^{\tau-t}_W|_{C^\sigma} &\le& C(\sigma) \left( |\tilde R^t_\tau|_{C^\sigma}|\partial_q \phi_W^{\tau-t}|^\sigma_{C^0} + |\tilde R^t_\tau|_{C^1} |\partial_q \phi_W^{\tau-t}|_{C^{\sigma -1}}+ |R^t_\tau|_{C^0}\right)
\end{eqnarray*}
and replacing the latter into the above integral 
\begin{eqnarray*}
|\varkappa^t|_{C^\sigma} &\le& C(\sigma) \int_t^{+\infty} |R^t_\tau|_{C^0}|z^\tau|_{C^\sigma}|\partial_q \phi_W^{\tau-t}|^\sigma_{C^0} d\tau + C(\sigma)\int_t^{+\infty} |R^t_\tau|_{C^0}|z^\tau|_{C^\sigma}|\partial_q\phi_W^{\tau-t}|_{C^{\sigma-1}}d\tau \\
&+&C(\sigma)\int_t^{+\infty} |\tilde R^t_\tau|_{C^\sigma}|\partial_q \phi_W^{\tau-t}|^\sigma_{C^0}|z^\tau|_{C^\sigma}d\tau + C(\sigma)\int_t^{+\infty} |\tilde R^t_\tau|_{C^1}|\partial_q \phi_W^{\tau-t}|_{C^{\sigma -1}}|z^\tau|_{C^\sigma}d\tau\\
&+& C(\sigma)\int_t^{+\infty} | R^t_\tau|_{C^0} |z^\tau|_{C^\sigma}d\tau.
\end{eqnarray*}
It remains to estimate each integral on the right-hand side of the latter. But, first, we note that for all $t \in J_\upsilon$ and $x<0$
\begin{equation*}
\int_t^{+\infty} e^{x \tau}(\tau - t)d\tau = {e^{xt} \over x^2}, \hspace{7mm} \int_t^{+\infty} e^{x \tau}(\tau - t)^2d\tau = -2{e^{xt} \over x^3}
\end{equation*}
where the latter is obtained by integrating by part. Now, thanks to Lemma \ref{psiNC}, Lemma \ref{LemmaRNC},~\eqref{HElambdaNC} and the latter
\begin{align*}
 &\int_t^{+\infty} |R^t_\tau|_{C^0}|z^\tau|_{C^\sigma}|\partial_q \phi_W^{\tau-t}|^\sigma_{C^0} d\tau\\%
 &\le C(\sigma) |z|^\upsilon_{\sigma, \lambda} \int_t^{+\infty} e^{\left(\bar c_1 \sigma + c_0^R\right)|\partial_q W|_{C^0}(\tau - t)} e^{-\lambda \tau} d\tau= C(\sigma){|z|^\upsilon_{\sigma, \lambda}  \over  \lambda - \left(\bar c_1 \sigma + c_0^R\right)|\partial_q W|_{C^0}}e^{\lambda t}\\
 &\int_t^{+\infty} |R^t_\tau|_{C^0}|z^\tau|_{C^\sigma}|\partial_q\phi_W^{\tau-t}|_{C^{\sigma-1}}d\tau\\%
 &\le C(\sigma) |z|^\upsilon_{\sigma, \lambda} \int_t^{+\infty} \left(1 + |\partial_q W|_{C^{\sigma-1}} (\tau - t)\right)e^{\left(c_\sigma + c_0^R\right)|\partial_q W|_{C^0}(\tau - t)} e^{-\lambda \tau} d\tau\\
 &=C(\sigma) |z|^\upsilon_{\sigma, \lambda} \int_t^{+\infty} e^{\left(c_\sigma + c_0^R\right)|\partial_q W|_{C^0}(\tau - t)} e^{-\lambda \tau} d\tau\\
 &+C(\sigma) |z|^\upsilon_{\sigma, \lambda} |\partial_q W|_{C^{\sigma-1}} \int_t^{+\infty} e^{\left(c_\sigma + c_0^R\right)|\partial_q W|_{C^0}(\tau - t)} e^{-\lambda \tau} (\tau - t)d\tau\\
 &=C(\sigma){|z|^\upsilon_{\sigma, \lambda}  \over  \lambda - \left(c_\sigma + c_0^R\right)|\partial_q W|_{C^0}}e^{\lambda t} + C(\sigma){|z|^\upsilon_{\sigma, \lambda}  |\partial_q W|_{C^{\sigma-1}}  \over \left( \lambda - \left(c_\sigma + c_0^R\right)|\partial_q W|_{C^0}\right)^2}e^{\lambda t}
 \end{align*}
\begin{align*}
&\int_t^{+\infty} |\tilde R^t_\tau|_{C^\sigma}|\partial_q \phi_W^{\tau-t}|^\sigma_{C^0}|z^\tau|_{C^\sigma}d\tau \\%
&\le C(\sigma) |z|^\upsilon_{\sigma, \lambda} \int_t^{+\infty} \left(1 + |\partial_q W|_{C^\sigma} (\tau - t)\right)e^{\left(c^R_\sigma + \bar c_1 \sigma\right)|\partial_q W|_{C^0}(\tau - t)} e^{-\lambda \tau} d\tau\\
&+ C(\sigma) |z|^\upsilon_{\sigma, \lambda} \int_t^{+\infty} |\partial_q W|_{C^1} |\partial_q W|_{C^{\sigma-1}} (\tau - t)^2 e^{\left(c^R_\sigma + \bar c_1 \sigma\right)|\partial_q W|_{C^0}(\tau - t)} e^{-\lambda \tau} d\tau\\
&=C(\sigma){|z|^\upsilon_{\sigma, \lambda}  \over  \lambda - \left(c^R_\sigma + \bar c_1 \sigma\right)|\partial_q W|_{C^0}}e^{\lambda t} + C(\sigma){|z|^\upsilon_{\sigma, \lambda}  |\partial_q W|_{C^{\sigma}}  \over \left( \lambda - \left(c^R_\sigma + \bar c_1 \sigma\right)|\partial_q W|_{C^0}\right)^2}e^{\lambda t}\\
&+C(\sigma){|z|^\upsilon_{\sigma, \lambda}  |\partial_q W|_{C^1} |\partial_q W|_{C^{\sigma-1}}  \over \left( \lambda - \left(c^R_\sigma + \bar c_1 \sigma\right)|\partial_q W|_{C^0}\right)^3}e^{\lambda t}
\end{align*}
\begin{align*}
&\int_t^{+\infty} |\tilde R^t_\tau|_{C^1}|\partial_q \phi_W^{\tau-t}|_{C^{\sigma -1}}|z^\tau|_{C^\sigma}d\tau\\%
&\le C(\sigma) |z|^\upsilon_{\sigma, \lambda} \int_t^{+\infty} \left(1 + |\partial_q W|_{C^1} (\tau - t)\right)\left(1 + |\partial_q W|_{C^{\sigma-1}} (\tau - t)\right)e^{\left(\bar c^R_1 + c_\sigma\right)|\partial_q W|_{C^0}(\tau - t)} e^{-\lambda \tau} d\tau\\
&\le C(\sigma) |z|^\upsilon_{\sigma, \lambda} \int_t^{+\infty} e^{\left(\bar c^R_1 + c_\sigma\right)|\partial_q W|_{C^0}(\tau - t)} e^{-\lambda \tau} d\tau\\
&+C(\sigma) |z|^\upsilon_{\sigma, \lambda} |\partial_q W|_{C^1} \int_t^{+\infty}(\tau - t)e^{\left(\bar c^R_1 + c_\sigma\right)|\partial_q W|_{C^0}(\tau - t)} e^{-\lambda \tau} d\tau\\
&+C(\sigma) |z|^\upsilon_{\sigma, \lambda}  |\partial_q W|_{C^{\sigma-1}} \int_t^{+\infty} (\tau - t)e^{\left(\bar c^R_1 + c_\sigma\right)|\partial_q W|_{C^0}(\tau - t)} e^{-\lambda \tau} d\tau\\
&+C(\sigma) |z|^\upsilon_{\sigma, \lambda} |\partial_q W|_{C^1}  |\partial_q W|_{C^{\sigma-1}} \int_t^{+\infty}  (\tau - t)^2 e^{\left(\bar c^R_1 + c_\sigma\right)|\partial_q W|_{C^0}(\tau - t)} e^{-\lambda \tau} d\tau\\
&=C(\sigma){|z|^\upsilon_{\sigma, \lambda}  \over  \lambda - \left(\bar c^R_1 + c_\sigma\right)|\partial_q W|_{C^0}}e^{\lambda t} + C(\sigma){|z|^\upsilon_{\sigma, \lambda} \left(|\partial_q W|_{C^1}+ |\partial_q W|_{C^{\sigma-1}} \right)  \over \left( \lambda - \left(\bar c^R_1 + c_\sigma\right)|\partial_q W|_{C^0}\right)^2}e^{\lambda t}\\
&+ C(\sigma){|z|^\upsilon_{\sigma, \lambda}  |\partial_q W|_{C^1} |\partial_q W|_{C^{\sigma-1}}  \over \left( \lambda - \left(\bar c^R_1 + c_\sigma\right)|\partial_q W|_{C^0}\right)^3}e^{\lambda t}\\
\end{align*}
\begin{align*}
&\int_t^{+\infty} | R^t_\tau|_{C^0} |z^\tau|_{C^\sigma}d\tau\\%
&\le C|z|^\upsilon_{\sigma, \lambda} \int_t^{+\infty}e^{c^R_0 |\partial_q W|_{C^0}(\tau - t)} e^{-\lambda \tau} d\tau = C(\sigma){|z|^\upsilon_{\sigma, \lambda}  \over \lambda -  c^R_0 |\partial_q W|_{C^0}}e^{\lambda t}.
\end{align*}
Now, we recall that $c_\sigma^\varkappa = \max\{c_\sigma + c_0^R, c_\sigma^R + c_0^R, \bar c_1^R + c_\sigma\}$. Hence, thanks to the latter
\begin{eqnarray*}
|\varkappa^t|_{C^\sigma} e^{\lambda t}&\le&  C(\sigma)\Bigg({1 \over  \lambda - \left(\bar c_1 \sigma + c_0^R\right)|\partial_q W|_{C^0}} + {1 \over  \lambda - \left(c_\sigma + c_0^R\right)|\partial_q W|_{C^0}}\\
&+& {|\partial_q W|_{C^{\sigma-1}}  \over \left( \lambda - \left(c_\sigma + c_0^R\right)|\partial_q W|_{C^0}\right)^2} + {1 \over  \lambda - \left(c^R_\sigma + \bar c_1 \sigma\right)|\partial_q W|_{C^0}}\\
&+&{|\partial_q W|_{C^{\sigma}}  \over \left( \lambda - \left(c^R_\sigma + \bar c_1 \sigma\right)|\partial_q W|_{C^0}\right)^2} + {|\partial_q W|_{C^1} |\partial_q W|_{C^{\sigma-1}}  \over \left( \lambda - \left(c^R_\sigma + \bar c_1 \sigma\right)|\partial_q W|_{C^0}\right)^3}\\
&+& {1 \over  \lambda - \left(\bar c^R_1 + c_\sigma\right)|\partial_q W|_{C^0}} + { |\partial_q W|_{C^1}+ |\partial_q W|_{C^{\sigma-1}}  \over \left( \lambda - \left(\bar c^R_1 + c_\sigma\right)|\partial_q W|_{C^0}\right)^2}\\
&+& {|\partial_q W|_{C^1} |\partial_q W|_{C^{\sigma-1}}  \over \left( \lambda - \left(\bar c^R_1 + c_\sigma\right)|\partial_q W|_{C^0}\right)^3} + {1  \over \lambda -  c^R_0 |\partial_q W|_{C^0}}\Bigg)|z|^\upsilon_{\sigma, \lambda}\\
&\le& C(\sigma) \Bigg({1 \over  \lambda - c_\sigma^\varkappa|\partial_q W|_{C^0}} + {|\partial_q W|_{C^\sigma}  \over \left( \lambda - c_\sigma^\varkappa|\partial_q W|_{C^0}\right)^2} + {|\partial_q W|_{C^1} |\partial_q W|_{C^{\sigma-1}}  \over \left( \lambda - c^\varkappa_\sigma|\partial_q W|_{C^0}\right)^3} \Bigg)|z|^\upsilon_{\sigma, \lambda}
\end{eqnarray*}
for all $t \in J_\upsilon$. Furthermore, we conclude the proof by taking the sup for all $t \in J_\upsilon$ on the left-hand side of the latter. 
\end{proof}

We observe that we do not find a unique solution to the homological equation, unlike when $W$ is constant. That is why we prove only the existence of a right inverse of the differential of $\mathcal{F}$ introduced in Section \ref{OTPNC}. 

Concerning the regularity, we are not able to find holomorphic solutions to the previous homological equation. Here, we solve~\eqref{HENC} thanks to a suitable change of coordinates~\eqref{hNC} that rectifies the dynamics on the torus. Unfortunately, this change of coordinates depends on the flow of $W$, and Lemma \ref{psiNC} shows that the derivatives of this flow diverge exponentially fast in time. This prevents us from well defining the change of coordinate~\eqref{hNC} in the case of holomorphic functions defined on complex neighbourhoods of the phase space.

Also, in the case of $C^\infty$ functions, we cannot find $C^\infty$ solutions of~\eqref{HENC}. Because the bigger $\sigma$ is, the more we have to take $\lambda$ big in order to solve the homological equation. 

\subsection{Regularity of $\mathcal{F}$}\label{RFNC}

We recall the definition of the functional $\mathcal{F}$, which we properly introduced in Section \ref{PSNC}.
\begin{equation*}
\mathcal{F} :  \mathcal{A} \times \mathcal{B} \times \mathcal{M} \times \mathcal{M} \times \mathcal{W} \times  \mathcal{U} \times  \mathcal{V} \longrightarrow \mathcal{Z} \times \mathcal{G}\\
\end{equation*}
\begin{equation*}
\mathcal{F}(a, b,m , \bar m, W, u, v) = (F_1(b, \bar m, u,v), F_2(a,b,m,W,u,v))
\end{equation*}
with 
\begin{eqnarray*}
F_1(b, \bar m,W,u,v) &=& W \circ  (\mathrm{id} + u) - W + b \circ \tilde u + \left(\bar m \circ \tilde \varphi  \right) v - \left(\nabla u \right)\bar W \\,
&=& \left(\int_0^1\partial_q W(\mathrm{id} + \tau u) d\tau \right)u+ b\circ \tilde u +  \left(\bar m \circ \tilde \varphi  \right) v - \left(\nabla u \right)\bar W \\
F_2(a,b, m,W,u,v) &=& \partial_q a\circ \tilde u  + \left(\partial_q W \circ (\mathrm{id} + u) + \partial_q b \circ \tilde u \right) v\\ 
&+& \left(\partial_q m \circ \tilde  \varphi   \right) \cdot v^2 + \left(\nabla v \right)\bar W,
\end{eqnarray*}
where the second line of the latter is a consequence of the Taylor formula. Thanks to Proposition \ref{normpropertiesNC}, the functional $\mathcal{F}$ is well-defined and continuous. Moreover,  $\mathcal{F}$ is differentiable with respect to the components $(u,v)$, with
\begin{eqnarray*}
D_{(u,v)} F_1(b,\bar m, W,u,v)(\hat u, \hat  v) &=& D_u F_1(b,\bar m, W,u,v)\hat u  + D_v F_1(b,\bar m, W,u,v)\hat v  \\
&=& \left(\partial_q W \circ(\mathrm{id} + u) + \partial_q b \circ \tilde u \right) \hat u + v^T \left(\partial_q \bar m \circ \tilde  \varphi \right)\hat u\\
&+& v^T \left( \partial_p \bar m \circ \tilde  \varphi \right)\hat v + \left(\bar m \circ \tilde  \varphi \right)\hat v - \left(\nabla \hat u \right) \bar W\\
D_{(u,v)} F_2(a, b,m,W,u,v) (\hat u, \hat v) &=& D_u F_2(a, b,m,W,u,v) \hat u  + D_v F_2(a, b,m,W,u,v) \hat v   \\
&=& \left(\partial^2_q a \circ \tilde u \right) \hat u  + v^T \left( \partial^2_q W\circ (\mathrm{id}+u) + \partial^2_q b \circ \tilde u \right) \hat u\\
&+&  (v^T)^2 \left( \partial^2_q m \circ \tilde  \varphi \right)\hat u + \left(\partial_q W \circ (\mathrm{id}+u)+ \partial_q b \circ \tilde u \right) \hat v\\
&+&  (v^T)^2 \left(\partial^2_{pq} m \circ \tilde  \varphi \right)\hat v + 2 v^T \left(\partial_q m \circ \tilde  \varphi \right)\hat v\\
&+&  \left(\nabla \hat v \right)\bar W,
\end{eqnarray*}
where $T$ stands for transpose. These differentials are continuous. Furthermore, the differential $D_{(u,v)}\mathcal{F}$ calculated in $(0,0,m,\bar m, W, 0, 0)$ is equal to 
\begin{equation}
\label{DuvFNC}
D_{(u,v)} \mathcal{F}(0,0,m, \bar m,W, 0,0) (\hat u, \hat v) = \begin{pmatrix} \partial_q W\hat u - \left(\nabla \hat u\right)\bar W + \bar m_0\hat v\\
\partial_q W\hat v +\left(\nabla \hat v\right)\bar W \end{pmatrix},
\end{equation}
where for all $(q, t) \in \T^n \times J_{\upsilon'}$ we let $\bar m_0(q,t) = \bar m(q,0,t)$. The following lemma proves that, for all fixed $m$, $\bar m \in \mathcal{M}$ and $W \in \mathcal{W}$, $D_{(u,v)} \mathcal{F}(0,0,m, \bar m,W, 0,0)$ admits a right inverse.

\begin{lemma}
\label{lemminvNC}
For all $(z, g) \in \mathcal{Z} \times \mathcal{G}$, there exists $(\hat u, \hat v) \in \mathcal{U} \times \mathcal{V}$ such that 
\begin{equation}
\label{DuvFLemmaNC}
D_{(u,v)} \mathcal{F}(0,0,m, \bar m,W, 0,0) (\hat u, \hat v)= (z,g).
\end{equation}
Moreover, for a suitable constant $\bar C$ depending on $\sigma$, $\Upsilon$, $\lambda$, $|\partial_q W|_{C^0}$ and $|\partial_q W|_{C^\sigma}$
\begin{equation*}
|\hat u| \le \bar C\left(|g|^{\upsilon'}_{\sigma, \lambda} + |z|^{\upsilon'}_{\sigma, \lambda}\right), \quad |\hat v| \le \bar C|g|^{\upsilon'}_{\sigma, \lambda}
\end{equation*}
where, we recall that $|\hat u| = \max\{|\hat  u|^{\upsilon'}_{\sigma, \lambda}, | \left(\nabla \hat  u \right) \bar W|^{\upsilon'}_{\sigma, \lambda} \}$ and \newline $|\hat v| = \max\{|\hat v|^{\upsilon'}_{\sigma, \lambda}, | \left(\nabla \hat v \right) \bar W|^{\upsilon'}_{\sigma, \lambda} \}$.
\end{lemma}
\begin{proof}
The proof of this lemma relies on~\eqref{lambdaNC} and Lemma \ref{homoeqlemmaNC}. Thanks to~\eqref{DuvFNC}, we can rewrite equation~\eqref{DuvFLemmaNC} in terms of the following system in the unknown $(\hat u, \hat v)$
\begin{equation}
\label{DuvFsystemNC}
\begin{cases}
\left(\nabla \hat u\right)\bar W  -  \partial_q W\hat u = \bar m_0\hat v - z\\
\left(\nabla \hat v\right)\bar W + \partial_q W\hat v = g.
\end{cases}
\end{equation}
These equations are decoupled, and hence we can study them separately. We begin by solving the last one. Then, we replace the found solution $\hat v$ in the first equation, which now can be solved, and we conclude the proof of this lemma.

By~\eqref{lambdaNC} and Lemma \ref{homoeqlemmaNC}, a solution $\hat v$ of the second equation of the above system exists and satisfies 
\begin{equation*}
|\hat v|^{\upsilon'}_{\sigma, \lambda} \le C(\sigma, \lambda, |\partial_q W|_{C^0}, |\partial_q W|_{C^\sigma})|g|^{\upsilon'}_{\sigma, \lambda}.
\end{equation*}
 Moreover, thanks to Proposition \ref{normpropertiesNC},~\eqref{DuvFsystemNC} and the latter
\begin{eqnarray*}
| \left(\nabla \hat v \right) \bar W|^{\upsilon'}_{\sigma, \lambda} &=& |g - \partial_q W\hat v |^{\upsilon'}_{\sigma, \lambda} \le |g|^{\upsilon'}_{\sigma, \lambda} + C(\sigma)|\partial_q W|_{C^\sigma}|\hat v |^{\upsilon'}_{\sigma, \lambda}\\
&\le&C(\sigma, \lambda, |\partial_q W|_{C^0}, |\partial_q W|_{C^\sigma})|g|^{\upsilon'}_{\sigma, \lambda}.
\end{eqnarray*}
As a consequence of the previous estimates, we obtain
\begin{equation*}
|\hat v| = \max\{|\hat v|^{\upsilon'}_{\sigma, \lambda}, | \left(\nabla \hat v \right) \bar W|^{\upsilon'}_{\sigma, \lambda} \} \le C(\sigma, \lambda, |\partial_q W|_{C^0}, |\partial_q W|_{C^\sigma})|g|^{\upsilon'}_{\sigma, \lambda},
\end{equation*}
which proves the first estimate of this lemma. Now, we can solve the first equation of~\eqref{DuvFsystemNC} where $\hat v$ is known. Thanks to Proposition \ref{normpropertiesNC} and the previous estimate
\begin{eqnarray*}
|\bar m_0\hat v - z|^{\upsilon'}_{\sigma, \lambda} &\le& C(\sigma) \Upsilon |\hat v|^{\upsilon'}_{\sigma, \lambda} + |z|^{\upsilon'}_{\sigma, \lambda}\le C(\sigma, \lambda, |\partial_q W|_{C^0}, |\partial_q W|_{C^\sigma}) \Upsilon |g|^{\upsilon'}_{\sigma, \lambda} + |z|^{\upsilon'}_{\sigma, \lambda}\\
&\le& C(\sigma, \lambda, \Upsilon, |\partial_q W|_{C^0}, |\partial_q W|_{C^\sigma})\left(|g|^{\upsilon'}_{\sigma, \lambda} + |z|^{\upsilon'}_{\sigma, \lambda}\right).
\end{eqnarray*}
By~\eqref{lambdaNC} and Lemma \ref{homoeqlemmaNC}, a solution of the first equation of~\eqref{DuvFsystemNC} exists and we have that
\begin{eqnarray*}
|\hat u|^{\upsilon'}_{\sigma, \lambda} &\le&C(\sigma, \lambda, |\partial_q W|_{C^0}, |\partial_q W|_{C^\sigma})|\bar m_0\hat v - z|^{\upsilon'}_{\sigma, \lambda} \\
&\le& C(\sigma, \lambda, \Upsilon, |\partial_q W|_{C^0}, |\partial_q W|_{C^\sigma})\left(|g|^{\upsilon'}_{\sigma, \lambda} + |z|^{\upsilon'}_{\sigma, \lambda}\right).
\end{eqnarray*}
 It remains to estimate $| \left(\nabla \hat  u \right) \bar W|^{\upsilon'}_{\sigma, \lambda}$ to conclude the proof of this lemma. This is a consequence of the previous estimates and~\eqref{DuvFsystemNC}
\begin{eqnarray*}
| \left(\nabla \hat u \right) \bar W|^{\upsilon'}_{\sigma, \lambda} &=& |\bar m_0\hat v - z +  \partial_q W\hat u|^{\upsilon'}_{\sigma, \lambda} \le |\bar m_0\hat v - z|^{\upsilon'}_{\sigma, \lambda} + C(\sigma)| \partial_q W|_{C^\sigma}|\hat u|^{\upsilon'}_{\sigma, \lambda}\\
&\le& C(\sigma, \lambda, \Upsilon, |\partial_q W|_{C^0}, |\partial_q W|_{C^\sigma})\left(|g|^{\upsilon'}_{\sigma, \lambda} + |z|^{\upsilon'}_{\sigma, \lambda}\right)
\end{eqnarray*}
and thus 
\begin{equation*}
|\hat u| = \max\{|\hat  u|^{\upsilon'}_{\sigma, \lambda}, | \left(\nabla \hat  u \right) \bar W|^{\upsilon'}_{\sigma, \lambda} \} \le C(\sigma, \lambda, \Upsilon, |\partial_q W|_{C^0}, |\partial_q W|_{C^\sigma})\left(|g|^{\upsilon'}_{\sigma, \lambda} + |z|^{\upsilon'}_{\sigma, \lambda}\right).
\end{equation*}
\end{proof}

\subsection{Existence of a $C^\sigma$-asymptotic torus}

This part is devoted to proving the existence of a $C^\sigma$-asymptotic torus associated to $(X_H, X_{\tilde h}, \varphi_0, W)$. To this end, we fix $x = (a,b)$, where $a$ and $b$ are those defined by~\eqref{H2NC}. It is straightforward to verify that $(a,b) \in \mathcal{A} \times \mathcal{B}$. Moreover, we define the following Banach space $\left(\mathcal{Y}, |\cdot|\right)$ such that $\mathcal{Y} = \mathcal{U} \times \mathcal{V}$ and, for all $y = (u,v) \in \mathcal{Y}$, $| y | = \max \{|u|, |v| \}$. Let $m$, $\bar m \in \mathcal{M}$ and $W\in\mathcal{W}$ be as in~\eqref{H2NC}, we rewrite $\mathcal{F}$ in the following form
\begin{eqnarray}
\label{FICNC}
\mathcal{F}(x,m, \bar m,W,y) &=& D_{(u,v)} \mathcal{F}(0,0,m, \bar m,W,0,0) y  + \mathcal{R}(x, m, \bar m, y).
\end{eqnarray}
For fixed $x$, $m$, $\bar m$ and $W$, the purpose of this section is to find $y \in \mathcal{Y}$ such that 
\begin{equation*}
\mathcal{F}(x,m, \bar m,W,y) =0. 
\end{equation*}
Let $\eta (m, \bar m,W)$ be the right inverse of $D_{(u,v)} \mathcal{F}(0,0,m, \bar m,W,0,0)$ whose existence is guaranteed by Lemma \ref{lemminvNC}. Therefore, we are looking for $y \in \mathcal{Y}$ in such a way that 
\begin{equation*}
y = y - \eta (m, \bar m,W) \mathcal{F}(x,m, \bar m,W,y).
\end{equation*}
To this end, we define the following functional 
\begin{equation*}
\mathcal{L}(x,m, \bar m,W, \cdot) :  \mathcal{Y} \longrightarrow  \mathcal{Y}
\end{equation*}
where
\begin{equation}
\label{LNC}
\mathcal{L}(x,m, \bar m, W,y) =  y - \eta (m, \bar m,W) \mathcal{F}(x,m, \bar m, W,y).
\tag{$\mathcal{L}$}
\end{equation}
It is well defined, and by the regularity of $\mathcal{F}$, we deduce that $\mathcal{L}$ is continuous and differentiable with respect to $y = (u,v)$ with differential $D_y\mathcal{L}$ continuous. The proof is reduced to find a fixed point of the latter. 
The following lemma is the main tool to conclude the proof of Theorem \ref{Thm1NC}. 

\begin{lemma}
\label{lemmautilethmNC}
There exists $\upsilon'$ large enough with respect to $n$, $\sigma$, $\lambda$, $|\partial_q W|_{C^{\sigma+1}}$ and $\Upsilon$, such that, for all  $y_*$,$y \in \mathcal{Y}$ with $|y_*|\le 1$,
\begin{equation}
\label{L2NC}
|D_y\mathcal{L}(x,m, \bar m,W,y_*) y| \le {1\over 2} |y|.
\end{equation}
\end{lemma}
\begin{proof}
The proof rests on Lemma \ref{lemminvNC}. By~\eqref{LNC}, for all $y$, $y^* \in \mathcal{Y}$
\begin{equation*}
D_y\mathcal{L}(x,m, \bar m,W,y_*) y= \left(\mathrm{Id} - \eta (m, \bar m,W)D_{(u,v)}\mathcal{F}(x,m, \bar m,W, y_*)\right)y.
\end{equation*}
We can reformulate the problem in terms of estimating the solution $\hat y = (\hat u, \hat v) \in \mathcal{Y}$ of the following system
\begin{align}
\label{systemfinalNC}
&D_{(u,v)} \mathcal{F}(0,0,m,\bar m,W,0,0)\hat y \\
&= \Big(D_{(u,v)}\mathcal{F}(0,0,m,\bar m,W,0,0) - D_{(u,v)}\mathcal{F}(x,m,\bar m,W, y_*)\Big) y. \nonumber
\end{align}
It suffices to estimate the right-hand side of the latter. Then, we conclude the proof by  Lemma \ref{lemminvNC}. We point out that $y_* = (u_*, v_*) \in \mathcal{Y}$ and for all $(q,t) \in \T^n \times J_{\upsilon'}$, we let
\begin{equation*}
\tilde u_*(q,t) = (q + u_*(q,t),t), \quad \tilde \varphi_*(q,t) = (q + u_*(q,t), v_*(q,t) ,t).
\end{equation*}
Thanks to~\eqref{DuvFNC}, we can rewrite the right-hand side of~\eqref{systemfinalNC} in the following form
\begin{equation*}
 \begin{pmatrix}\partial_q Wu - \left(\nabla u\right)\bar W + \bar m_0 v -  D_{(u,v)} F_1(b,\bar m, W,y_*)y\\ \partial_q Wv +\left(\nabla  v\right)\bar W -D_{(u,v)} F_2(x,m,W,y_*)y  \end{pmatrix}
\end{equation*}
(see Section \ref{RFNC}), moreover 
\begin{eqnarray*}
\partial_q Wu - \left(\nabla u\right)\bar W + \bar m_0 v -  D_{(u,v)} F_1(b,\bar m, W,y_*)y&=& \left( \bar m_0 - \bar m \circ \tilde  \varphi_* \right) v - \left(\partial_q b \circ \tilde u_* \right)  u \\
&+& \left(\partial_qW - \partial_q W\circ (\mathrm{id} + u_*)\right)u\\
&-& v_*^T \left(\partial_q \bar m \circ \tilde  \varphi_* \right) u - v_*^T \left( \partial_p \bar m \circ \tilde  \varphi_* \right) v \\
\partial_q Wv +\left(\nabla  v\right)\bar W -D_{(u,v)} F_2(x,m,W,y_*)y  &=& \left(\partial_q W - \partial_q W \circ \left(\mathrm{id} + u_*\right)\right)v\\
&-&\left(\partial^2_q a \circ \tilde u_* \right)  u \\
&-& v_*^T \left( \partial^2_q W \circ (\mathrm{id} + u_*) + \partial^2_q b \circ \tilde u_* \right) u \\
&-&  (v_*^T)^2 \left( \partial^2_q m \circ \tilde  \varphi_* \right) u - \left(\partial_q b \circ \tilde u_* \right) v \\
&-&  (v_*^T)^2 \left(\partial^2_{pq} m \circ \tilde  \varphi_* \right)v \\
&-& 2 v_*^T \left(\partial_q m \circ \tilde  \varphi_* \right)v.
\end{eqnarray*}
 Now, thanks to property \textit{2.} of Proposition \ref{Holder2}, 
\begin{eqnarray*}
 \left|\left(\partial_q Wu - \left(\nabla u\right)\bar W + \bar m_0 v -  D_{(u,v)} F_1(b,\bar m, W,y_*)y\right)^t\right|_{C^\sigma} &\le& C(\sigma)\Big( \left|\left(\bar m_0^t -\bar m \circ \tilde  \varphi_* \right)^t\right|_{C^\sigma}  \left|v^t\right|_{C^\sigma} \\
&+& \left|\left(\partial_q b \circ \tilde u_* \right)^t\right|_{C^\sigma}   \left|u^t \right|_{C^\sigma}  \\
&+&  |\left(\partial_qW - \partial_q W\circ (\mathrm{id} + u_*)\right)^t|_{C^\sigma} |u^t|_{C^\sigma} \\
&+& |v^t_*|_{C^\sigma}  \left|\left(\partial_q \bar m \circ \tilde  \varphi_* \right)^t\right|_{C^\sigma}  |u^t|_{C^\sigma}   \\
&+&  |v^t_*|_{C^\sigma}  \left| \left(\partial_p \bar m \circ \tilde  \varphi_* \right)^t\right|_{C^\sigma}   |v^t|_{C^\sigma}  \Big)
\end{eqnarray*}
for all $t \in J_{\upsilon'}$. We have to estimate each term on the right-hand side of the latter. We begin with the third one. 
\begin{align*}
&|\left(\partial_qW - \partial_q W\circ (\mathrm{id} + u_*)\right)^t|_{C^\sigma} |u^t|_{C^\sigma} \le C(\sigma)| \partial^2_q W\circ (\mathrm{id} + \tau u_*)^tu^t_*|_{C^\sigma} |u^t|_{C^\sigma} \\
&\le C(\sigma)|\partial_q W|_{C^{\sigma+1}}|y_*| e^{-\lambda t}\left(1 + \left(1 + |\partial_q u_*^t|_{C^0}\right)^\sigma + |\partial_q u_*^t|_{C^{\sigma-1}}\right)|y|e^{-\lambda t} \\
&\le C(\sigma)|\partial_q W|_{C^{\sigma+1}}e^{-\lambda \upsilon'} |y|e^{-\lambda t}
\end{align*}
for all $t \in J_{\upsilon'}$. The first line of the latter is due to the mean value theorem for a suitable $\tau \in [0,1]$. In the second line, we use properties \textit{2.} and \textit{5.} of Proposition \ref{Holder2}. The last line is due to $|y_*| \le 1$. Similarly to the previous case, thanks to the mean value theorem, properties \textit{2.} and \textit{5.} of Proposition \ref{Holder2} and $|y_*| \le 1$, we obtain 
\begin{eqnarray*}
\left|\left(\bar m_0 -\bar m \circ \tilde  \varphi_* \right)^t\right|_{C^\sigma} |v^t|_{C^\sigma} &\le& C(\sigma) \Big(|\partial_q \bar m^t(\mathrm{id} + \tau u_*, \tau v_*)u^t_*|_{C^\sigma}\\
&+& |\partial_p \bar m^t(\mathrm{id} + \tau u_*, \tau v_*)v^t_*|_{C^\sigma}\Big) |v^t|_{C^\sigma} \\
&\le&C(\sigma)\Upsilon \left(|u^t_*|_{C^\sigma} + |v^t_*|_{C^\sigma}\right)|v^t|_{C^\sigma}\\
&\le&C(\sigma)\Upsilon e^{-\lambda \upsilon'} |y|e^{-\lambda t}\\
\left|\left(\partial_q b \circ \tilde u_* \right)^t\right|_{C^\sigma}  \left|u^t \right|_{C^\sigma} &\le& C(\sigma) |b|^{\upsilon^*}_{\sigma+2,\lambda} e^{-\lambda \upsilon'}|y|e^{-\lambda t}\\
&\le& C(\sigma) \Upsilon e^{-\lambda t}|y|e^{-\lambda t}\\
 |v_*^t|_{C^\sigma} \left|\left(\partial_q \bar m \circ \tilde  \varphi_* \right)^t\right|_{C^\sigma} |u^t|_{C^\sigma} &\le&C(\sigma) |y_*|e^{-\lambda \upsilon'}\Upsilon |y|e^{-\lambda t} \\
 |v_*^t|_{C^\sigma} \left|\left( \partial_p \bar m \circ \tilde  \varphi_* \right)^t\right|_{C^\sigma} |v^t|_{C^\sigma} &\le&C(\sigma) |y_*|e^{-\lambda \upsilon'}\Upsilon |y|e^{-\lambda t} 
\end{eqnarray*}
for all $t \in J_{\upsilon'}$. Therefore, for $\upsilon'$ large enough, the above estimates imply
\begin{eqnarray*}
 \left|\left(\partial_q Wu - \left(\nabla u\right)\bar W + \bar m_0 v -  D_{(u,v)} F_1(b,\bar m, W,y_*)y\right)^t\right|_{C^\sigma} &\le& {1 \over 4 \bar C}|y|e^{-\lambda t}
\end{eqnarray*}
for all $t \in J_{\upsilon'}$. We point out that $\bar C$ is the constant introduced in Lemma \ref{lemminvNC}.  Multiplying both sides of the latter by $e^{\lambda t}$ and taking the sup for all $t \in J_{\upsilon'}$, we obtain
\begin{eqnarray*}
 \left|\partial_q Wu - \left(\nabla u\right)\bar W + \bar m_0 v -  D_{(u,v)} F_1(b,\bar m, W,y_*)y\right|^{\upsilon'}_{\sigma, \lambda} &\le& {1 \over 4 \bar C}|y|.
\end{eqnarray*}
Similarly to the previous case, for $\upsilon'$ large enough, we have 
\begin{eqnarray*}
 \left|\partial_q Wv +\left(\nabla  v\right)\bar W -D_{(u,v)} F_2(x,m,W,y_*)y \right|^{\upsilon'}_{\sigma, \lambda} &\le&  {1 \over 4\bar C}|y|.
\end{eqnarray*}
This concludes the proof of this lemma. Thanks to Lemma \ref{lemminvNC}, a solution $\hat y \in \mathcal{Y}$ of~\eqref{systemfinalNC} exists satisfying
\begin{eqnarray*}
|\hat u| &\le& \bar C\Big(  \left|\partial_q Wu - \left(\nabla u\right)\bar W + \bar m_0 v -  D_{(u,v)} F_1(b,\bar m, W,y_*)y\right|^{\upsilon'}_{\sigma, \lambda}\\
&+& \left|\partial_q Wu - \left(\nabla u\right)\bar W + \bar m_0 v -  D_{(u,v)} F_1(b,\bar m, W,y_*)y\right|^{\upsilon'}_{\sigma, \lambda} \Big) \le {1 \over 2}|y|\\
|\hat v| &\le&\bar C \left|\partial_q Wu - \left(\nabla u\right)\bar W + \bar m_0 v -  D_{(u,v)} F_1(b,\bar m, W,y_*)y\right|^{\upsilon'}_{\sigma, \lambda}\le {1 \over 4}|y|
\end{eqnarray*}
and hence
\begin{equation*}
|D_y\mathcal{L}(x,m, \bar m,W,y_*) y| \le {1\over 2} |y|.
\end{equation*}
\end{proof}

We proved that $\mathcal{L}(x, m, \bar m, \cdot)$ is a contraction of a compact subset of $\mathcal{Y}$. So then, there exists a unique fixed point $y \in \mathcal{Y}$ with $|y| \le 1$. 

More specifically, there exists $(u, v) \in \mathcal{U} \times \mathcal{V}$ such that, for all $(q,t) \in \T^n \times J_{\upsilon'}$
\begin{equation*}
\varphi^t(q) = (q + u(q,t), v(q,t))
\end{equation*}
is a $C^\sigma$-asymptotic torus associated to $(X_H, X_{\tilde h}, \varphi_0, W)$. We conclude the proof by verifying that $\varphi^t$ is a Lagrangian $C^\sigma$-asymptotic torus.  Let $\psi_{t_0,H}^t$ and $\phi_{t_0,W}^t$ be the flows at time $t$ with initial time $t_0$ of $H$ and $W$, respectively. 

\begin{lemma}
$\varphi^{t_0}$ is Lagrangian for all $t_0 \in J_{\upsilon'}$.
\end{lemma}
\begin{proof}
Let $\alpha = dp \wedge dq$ be the standard symplectic form associated with $(q,p) \in \T^n \times B$. By~\eqref{hyp1bissNC}, we know that for all $t_0 \in J_{\upsilon'}$ and $t \ge 0$
\begin{equation}
\label{hyp1bissbissMD}
\psi_{t_0,H}^{t_0 + t} \circ \varphi^{t_0} = \varphi^{t_0 + t} \circ \phi_{t_0, W}^{t_0+t},
\end{equation}
and taking the pull-back with respect to the standard form $\alpha$ on both sides of the latter, we obtain
\begin{equation*}
(\varphi^{t_0})^*(\psi_{t_0,H}^{t_0 + t})^* \alpha =  ( \phi_{t_0, W}^{t_0+t})^*(\varphi^{t_0 + t})^* \alpha.
\end{equation*}
For all fixed $t$, $t_0 \in J_{\upsilon'}$, the flow $\psi^t_{t_0,H}$  is a symplectomorphisms. Then, $(\psi_{t_0, H}^t)^*\alpha = \alpha$ for all fixed $t$, $t_0 \in J_{\upsilon'}$. This implies that we can rewrite the latter as follows
\begin{equation*}
(\varphi^{t_0})^* \alpha =  ( \phi_{t_0, W}^{t_0+t})^*(\varphi^{t_0 + t})^* \alpha
\end{equation*}
for all fixed $t_0 \in J_{\upsilon'}$ and $t \ge 0$. We want to prove that, $\left((\varphi^{t_0})^* \alpha\right)_q = 0$ for all $q\in \T^n$. To this end, we observe that we can rewrite the right-hand side of the latter as follows
\begin{equation*}
\left((\phi_{t_0, W}^{t_0+t})^*(\varphi^{t_0 + t})^* \alpha\right)_q = \sum_{1 \le i <j \le n} \sum_{1 \le k < d \le n} \alpha^t_{i,j,k,d}(q) dq_k \wedge dq_d
\end{equation*}
for all $q \in \T^n$, where
\begin{eqnarray*}
\alpha^t_{i,j,k,d}(q)  &=& \left(\partial_{q_i}v^{t_0 + t} \cdot \partial_{q_j} \left(\mathrm{id} + u^t\right) - \partial_{q_j}v^{t_0 + t} \cdot \partial_{q_i} \left(\mathrm{id} + u^t\right)\right) \circ \phi_{t_0, W}^{t_0+t}(q) \\
&\times&\left(\partial_{q_k}\phi_{t_0, W, i}^{t_0+t}(q)\partial_{q_d}\phi_{t_0, W,j}^{t_0+t}(q) - \partial_{q_d}\phi_{t_0, W, i}^{t_0+t}(q)\partial_{q_k}\phi_{t_0, W,j}^{t_0+t}(q)\right),
\end{eqnarray*}
$\times$ stands for the usual multiplication in $\R$ and $\phi^{t_0+t}_{t_0, W} = (\phi^{t_0+t}_{t_0, W, 1},...,\phi^{t_0+t}_{t_0, W,n})$. Then, for fixed $1 \le i <j \le n$, $1 \le k < d \le n$, by Lemma \ref{psiNC}
\begin{eqnarray*}
\left|\alpha^t_{i,j,k,d}\right|_{C^0}  &\le& \left|\left(\partial_{q_i}v^{t_0 + t} \cdot \partial_{q_j} \left(\mathrm{Id} + u^{t_0 + t} \right) - \partial_{q_j}v^{t_0 + t} \cdot \partial_{q_i} \left(\mathrm{Id} + u^{t_0 + t} \right)\right) \circ \psi_{t_0, W}^{t_0+t}(q) \right|_{C^0} \\
&\times&\left|\partial_{q_k}\psi_{t_0, W, i}^{t_0+t}\partial_{q_d}\psi_{t_0, W,j}^{t_0+t} - \partial_{q_d}\psi_{t_0, W, i}^{t_0+t}\partial_{q_k}\psi_{t_0, W,j}^{t_0+t}\right|_{C^0}\\
&\le& \left|\partial_{q_i}v^{t_0 + t} \cdot \partial_{q_j} \left(\mathrm{Id} + u^{t_0 + t} \right) - \partial_{q_j}v^{t_0 + t} \cdot \partial_{q_i} \left(\mathrm{Id} + u^{t_0 + t} \right) \right|_{C^0} \\
&\times&\left(\left|\partial_{q_k}\psi_{t_0, W, i}^{t_0+t}\right|_{C^0} \left|\partial_{q_d}\psi_{t_0, W,j}^{t_0+t}\right|_{C^0} + \left|\partial_{q_d}\psi_{t_0, W, i}^{t_0+t}\right|_{C^0}\left|\partial_{q_k}\psi_{t_0, W,j}^{t_0+t}\right|_{C^0}\right)\\
&=& \left|\partial_{q_i}v_j^{t_0 + t} + \partial_{q_i}v^{t_0 + t} \cdot \partial_{q_j} u^{t_0 + t}  - \partial_{q_j}v_i^{t_0 + t} -\partial_{q_j}v^{t_0 + t} \cdot \partial_{q_i}u^{t_0 + t}  \right|_{C^0} \\
&\times&\left(\left|\partial_{q_k}\psi_{t_0, W, i}^{t_0+t}\right|_{C^0} \left|\partial_{q_d}\psi_{t_0, W,j}^{t_0+t}\right|_{C^0} + \left|\partial_{q_d}\psi_{t_0, W, i}^{t_0+t}\right|_{C^0}\left|\partial_{q_k}\psi_{t_0, W,j}^{t_0+t}\right|_{C^0}\right)\\
&\le& C\left|\partial_q v^{t_0 + t}\right|_{C^0}\left(1 + \left|\partial_q u^{t_0 + t}\right|_{C^0}\right) \left|\partial_q\psi_{t_0, W}^{t_0+t}\right|^2_{C^0}\\
&\le& C\left|v^{t_0 + t}\right|_{C^1}\left(1 + \left|u^{t_0 + t}\right|_{C^1}\right) \left|\partial_q\psi_{t_0, W}^{t_0+t}\right|^2_{C^0} \le C e^{-\lambda (t_0 + t)}e^{2\bar c_1 |\partial_q W|_{C^0}t}
\end{eqnarray*}
for a suitable constant $C\ge 1$. Thanks to~\eqref{lambdaNC}, taking the limit for $t \to +\infty$ on both sides of the latter, the term in the last line converges to zero. This concludes the proof of this lemma.
\end{proof}

\section{Proof of Corollary \ref{CorNC}}\label{ProofCor1NC}

The proof is quite similar to that of Theorem \ref{Thm1NC}.
We are looking for a $C^\sigma$-asymptotic torus $\psi^t$ associated to $(Z, W, \mathrm{Id}, W)$. More specifically, for given $Z$, we are searching for $\upsilon' \ge 0$ sufficiently large and a suitable function $u : \T^n \times J_{\upsilon'} \to \R^n$ such that 
\begin{equation*}
\psi(q,t) = q + u(q,t)
\end{equation*}
and in addition $\psi$ and $u$ satisfy 
\begin{align}
\label{hyp1NC2Z}
&  Z(\psi(q, t), t) -  \partial_q \psi(q, t) W(q) - \partial_t \psi(q, t) = 0,\\
\label{hyp2NC2Z}
&  \lim_{t \to +\infty}  |u^t|_{C^\sigma} = 0
\end{align}
for all $(q, t) \in \T^n \times J_{\upsilon'}$.

We begin by introducing a suitable functional $\mathcal{F}$ given by~\eqref{hyp1NC2Z}. To this end, we define
\begin{equation*}
\tilde \psi(q,t) = (q + u(q,t), t), 
\end{equation*}
for all $(q,t) \in \T^n \times J_{\upsilon'}$. The composition of $Z$ with $\tilde \psi$ is equal to
\begin{equation*}
Z \circ \tilde \psi (q, t) = W \circ \psi (q, t) + P \circ \tilde \psi (q, t)
\end{equation*}
for all $(q,t) \in \T^n \times J_{\upsilon'}$, moreover
\begin{equation*}
\partial_q \psi(q, t)W(q) + \partial_t \psi(q, t) = W(q) + \partial_q u(q, t)W(q) + \partial_t u (q, t)
\end{equation*}
for all $(q,t) \in \T^n \times J_{\upsilon'}$. Then,  we can rewrite~\eqref{hyp1NC2Z} as follows 
\begin{eqnarray}
\label{F.0GDY}
P\circ \tilde\psi + (W \circ \psi - W) - \left(\nabla u \right) \bar W = 0.
\end{eqnarray}
This is the sum of functions defined for all $(q, t) \in \T^n \times J_{\upsilon'}$ or $q \in \T^n $. As usual, we have omitted the arguments $(q, t)$ and $q$ in order to achieve a more elegant form. For the sake of clarity, we recall that, during the proof of Theorem \ref{Thm1NC}, we have introduced the following notation
\begin{equation*}
\left(\nabla u \right) \bar W = \left(\partial_q u\right) W + \partial_t u.
\end{equation*}

Before the introduction of the functional $\mathcal{F}$, let $\sigma \ge 1$ be the positive parameter defined in Corollary \ref{CorNC}. For a suitable positive parameter $\upsilon'\ge0$ that we will choose large enough in Lemma \ref{LemmasceltaCorNotCons}, we introduce the following Banach spaces $\left(\mathcal{P}, |\cdot |\right)$, $\left(\mathcal{U}, |\cdot |\right)$, $\left(\mathcal{Z}, |\cdot |\right)$ and $\left(\mathcal{W}, |\cdot |\right)$  
\vspace{5mm}
\begin{eqnarray*}
\mathcal{P} &=& \Big\{P : \T^n \times J_{\upsilon'}  \to \R^n  \hspace{1mm}| \hspace{1mm} P \in \mathcal{\bar S}^{\upsilon'}_{\sigma, 1} , \hspace{1mm} \mbox{and} \hspace{1mm} |P| =|P|^{\upsilon'}_{\sigma+1, \lambda} < \infty\Big\}\\
\mathcal{U} &=& \Big\{u :\T^n \times J_{\upsilon'}  \to \R^n  \hspace{1mm}| \hspace{1mm} u,\left( \nabla u\right) \Omega \in \mathcal{S}^{\upsilon'}_{\sigma} \\
&& \mbox{and} \hspace{1mm} |u| = \max\{|u|^{\upsilon'}_{\sigma, \lambda}, | \left(\nabla u \right) \Omega|^{\upsilon'}_{\sigma, \lambda} \} < \infty\Big\}\\
\mathcal{Z} &=& \Big\{z : \T^n \times J_{\upsilon'}  \to \R^n  \hspace{1mm}| \hspace{1mm} z \in \mathcal{S}^{ \upsilon'}_{\sigma} , \hspace{1mm} \mbox{and} \hspace{1mm} |z| =|z|^{\upsilon'}_{\sigma, \lambda} < \infty\Big\}\\
\mathcal{W} &=& \Big\{ W:\T^n \to \R^n \hspace{1mm}| \hspace{1mm} W\in C^{\sigma+1}(\T^n)\hspace{1mm} \mbox{and} \hspace{1mm} |W| = |W|_{C^{\sigma+1}}<\infty \Big\}.
\end{eqnarray*} 

Now, thanks to~\eqref{F.0GDY} and the previous Banach spaces, we have everything we need to introduce the functional $\mathcal{F}$. Let $\mathcal{F}$ be the following functional
\begin{equation*}
\mathcal{F} :  \mathcal{P} \times \mathcal{W} \times  \mathcal{U} \longrightarrow \mathcal{Z} \\
\end{equation*}
\begin{equation*}
\mathcal{F}(P,W,u) = P\circ \tilde\psi + (W \circ \psi - W) - \left(\nabla u \right) \bar W.
\end{equation*}
 We observe that for all $W \in \mathcal{W}$
\begin{equation*}
\mathcal{F}(0, W, 0)=0.
\end{equation*}
Therefore, we can reformulate our problem in the following form. We fix $W \in \mathcal{W}$ and for $P \in \mathcal{P}$ sufficiently close to $0$, we are looking for $u \in \mathcal{U}$ in such a way that $\mathcal{F}(P,W,u) = 0$. 

Concerning the associated linearized problem, the differential of $\mathcal{F}$ with respect to the variable $u$ calculated in $(0,W,0)$ is equal to
\begin{equation*}
D_u \mathcal{F}(0,W,0)\hat u = \partial_q W \hat u -\left(\nabla \hat u \right) \bar W.
\end{equation*}

The functional $\mathcal{F}$ is well defined, continuous, differentiable with respect to the coordinate $u$ with $D_u\mathcal{F}(P,W,u)$ continuous. Moreover, by Lemma \ref{homoeqlemmaNC},  for all fixed $W \in \mathcal{W}$, $D_u \mathcal{F}(0,W,0)$ admits a right inverse $\eta(W)$. Then, $\mathcal{F}$ satisfies the hypotheses of the implicit function theorem. 

We fix $P \in \mathcal{P}$ and $W \in \mathcal{W}$ as in Corollary \ref{CorNC} and we define the following functional
\begin{equation*}
\mathcal{L}(P,W,\cdot):\mathcal{U} \longrightarrow \mathcal{U}
\end{equation*}
 in such a way that 
\begin{equation*}
\mathcal{L}(P,W,u) =  u - \eta(W)\mathcal{F}(P,W,u).
\end{equation*}
The proof of Corollary \ref{CorNC} is reduced to find a fixed point of the latter. Similarly to the proof of Lemma \ref{lemmautilethmNC}, we have the following lemma
\begin{lemma}
\label{LemmasceltaCorNotCons}
There exists $\upsilon'$ large enough with respect to $n$, $\sigma$, $\lambda$ and $|\partial_q W|_{C^\sigma}$, such that, for all  $u_*$,$u \in \mathcal{U}$ with $|u_*|\le 1$,
\begin{equation*}
|D_u\mathcal{L}(P,W,u_*) u| \le {1\over 2} |u|.
\end{equation*}
\end{lemma}
This concludes the proof of Corollary \ref{CorNC}.

\section*{Acknowledgement}

\textit{These results are part of my PhD thesis, which I prepared at Université Paris-Dauphine. I want to thank my thesis advisors, Abed Bounemoura and Jacques Féjoz. Without their advice and support, this work would not exist.}

\textit{This project has received funding from the European Union’s Horizon 2020 research and innovation programme under the Marie Skłodowska-Curie grant agreement No 754362}. \includegraphics[scale=0.01]{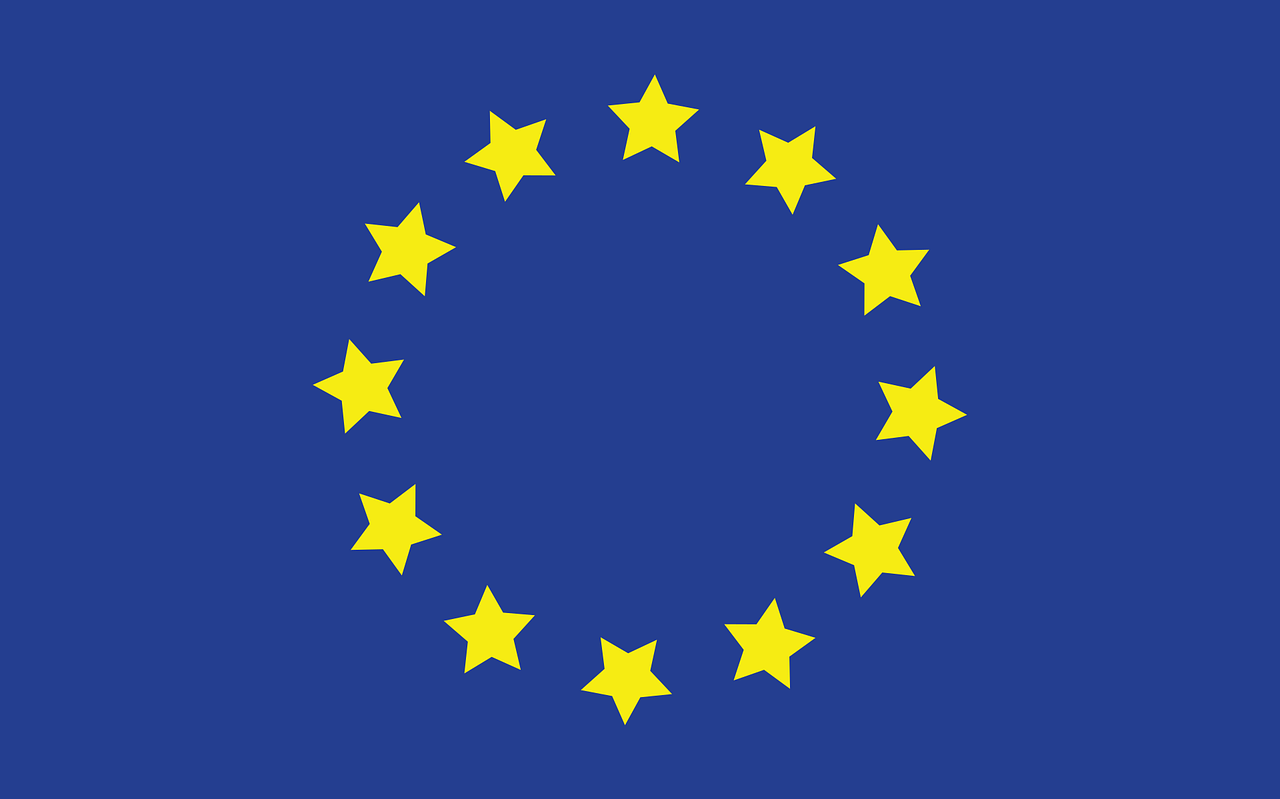}

\bibliographystyle{amsalpha}
\bibliography{Arbitraryfinal}

\end{document}